\documentclass[12pt]{article}

\newtheorem{df}{Definition}
\newtheorem{thm}[df]{Theorem}
\newtheorem{prop}[df]{Proposition}
\newtheorem{lem}[df]{Lemma}

\newtheorem{cor}[df]{Corollary}

\newcommand{\qed}{\hfill$\Box$}

\newcommand{\cir}{{\rm circ}}

\newcommand{\cbc}{{\rm cbc}}

\newcommand{\head}{{\frak h}}
\newcommand{\tail}{{\frak t}}

\usepackage{amssymb}

\usepackage{comment}
\usepackage{amsmath}
\usepackage[final]{graphicx}
\usepackage[dvipdfmx]{color}
\usepackage{tikz}

\begin{document}

%
%

\begin{center}
\LARGE{
	The Ihara expression for generalized weighted zeta functions of Bartholdi type 
	on finite digraphs
}

\vspace*{5mm}
\large{
	Ayaka Ishikawa\\
	Yokohama National University,\\
	Hodogaya, Yokohama 240-8501, Japan\\
	\vspace{3mm}
	Hideaki Morita\\
	Muroran Institute of Technology\\
	Mizumoto, Muroran 050-8585, Japan\\
	\vspace{3mm}
 	and\\
	Iwao Sato\\
	Oyama National College of Technology\\
	Nakakuki, Oyama 323-0806, Japan
	 }
\end{center}

\begin{abstract}
	The Ihara expression of a weighted zeta function 
	for a general finite digraph is given. 
	It unifies all the Ihara expressions obtained for known zeta functions for finite digraphs. 
	Any digraph in this paper permits multi-edges and multi-loops. 
\end{abstract}

	\section{Introduction}
	
	Graph zeta functions are formal power series associated with finite graphs, 
	that describe closed paths, cycles, or prime cycles of finite graphs. 
	The Ihara expression is one of the four expressions that graph zeta functions may have. 
	It is a determinant expression described by the adjacency matrix and the degree matrix of graphs, 
	and have been one of the main interest in the study of graph zeta functions 
	since its origin was discovered for the Ihara zeta function by Y. Ihara\cite{ihara66}. 
	Ihara's result  is given for the case where graphs are regular. 
	For a general finite graph, 
	the Ihara expression is obtained by H. Bass \cite{bass92}, 
	and the theorem is called the Bass-Ihara theorem. 
	Subsequently, the Bass-Ihara theorem is given its proof in various ways 
	(c.f., \cite{foatazeilberger99, kotanisunada00, northshield99, starkterras96}), 
	and also it is generalized to other graph zeta functions 
	(c.f., \cite{bartholdi99, mizunosato04, sato07, watanabefukumizu10}). 
	
	Recently, a new significance is given to the Ihara expression from an unexpected point of view. 
	In \cite{konnosato12}, N. Konno and I. Sato show that 
	the Ihara expression for the Sato zeta function \cite{sato07}, 
	a weighted version of the Ihara zeta function, gives a precise description 
	of the characteristic polynomial of the Grover matrix $U_{\Gamma}$ on a finite simple graph $\Gamma$. 
	The Grover matrix $U_{\Gamma}$ is the time-evolution matrix of the Grover walk \cite{grover96}, 
	the most extensively studied quantum walk model on finite graphs (c.f.,  \cite{AAKV01, HKSS17, kendon06}). 
	The spectrum of the time-evolution matrix is significant since it describes 
	fundamental features of the corresponding quantum walk: 
	mixing time, periodicity, and localization. 
	In \cite{EHSW06a}, Emms, Severini, Hancock and Wilson consider the Grover matrix on finite simple graphs, 
	and describe its spectrum directly considering the eigenvectors of $U_{\Gamma}$. 
	Although Emms et al.\cite{EHSW06a} do not present the characteristic polynomial explicitly, 
	Konno and Sato \cite{konnosato12} give its fine description, 
	provided by the viewpoint coming from the fact that $U_{\Gamma}$ is 
	realized as the edge matrix for the Sato zeta function of $\Gamma$. 
	This process corresdonds to reformulating the Hashimoto expression, 
	another determinant expression for graph zeta functions, 
	into the Ihara expression. 
	In other words, the Konno-Sato theorem implies that 
	the Hashimoto expression describes the time-evolution of the Grover walk 
	and the Ihara expression the spectral mapping theorem (c.f., \cite{matueogurisusegawa17, segawasuzuki19}). 
	Thus, the Ihara expression of graph zetas progresses its significance 
	from a quantum walk point of view.

	Though the construction of the Ihara expression has relied on case-by-case arguments 
	(c.f., \cite{bartholdi99, bass92, mizunosato04, sato07}), 
	Y. Watanabe and K. Fukumizu \cite{watanabefukumizu10} provided a new point of view 
	which has a possibility to unify previous studies. 
	Their idea is reformulating the Hashimoto expression into the Ihara expression 
	by linear algebraic method, 
	which essentially coincide with the idea of the Konno-Sato theorem. 
	Founded on this idea, 
	we construct the Ihara expression for the generalized weighted zeta function 
	by reformulating the Hashimoto expression. 
	The generalized weighted zeta function is defined in \cite{morita20} 
	which unifies the graph zeta functions that appeared in previous studies. 
	Also in \cite{morita20}, 
	it is investigated the condition for the Hashimoto expression to exist 
	and one can see that the condition is satisfied by the generalized weighted zeta function. 
	Thus, the generalized weighted zeta function has the Hashimoto expression, 
	and we can consider the problem to reformulate it into the Ihara expression. 
	The main result gives the Ihara expression for the generalized weighted zeta function 
	which is described by the weighted adjacency matrix and the weighted backtrack matrix for finite digraphs. 
	The weighted adjacency matrix and the backtrack matrix generalize 
	the adjacency matrix and the degree matrix respectively. 
	These generalized matrices coincides with the ordinary ones on the symmetric digraphs of a finite simple graph. 
	We will make one more remark on digraphs here. 
	In our development, 
	the definition of {\lq\lq inverse arcs\rq\rq} is significant. 
	Since our underlying digraph $\Delta$ may have multi-arcs and multi-loops, 
	there will be various definitions for inverse arcs. 
	In this article, every arc with the inverse direction to an arc $a$ is defined to be an inverse arc of $a$. 
	This definition of inverse arcs is a natural generalization 
	of the usual definition which we consider in the case where 
	$\Delta$ is the symmetric digraph $\Delta(\Gamma)$ for a finite simple graph $\Gamma$.

	The remaining part of this paper is organized as follows. 
	In Section 2, 
	we recall fundamental notation on graphs, digraphs, words, and dynamical systems, 
	which is required to define the combinatorial zeta functions. 
	In Section 3, 
	we review a fundamental theory of combinatorial zeta functions following \cite{morita20}. 
	We introduce the definition of a graph zeta function in general, 
	and see that it is combinatorial. 
	As a consequence, a graph zeta function has a determinant expression, 
	called the Hashimoto expression. 
	In Section 4, 
	we review the definition of the generalized weighted zeta function for finite digraphs, 
	and see that classical graph zeta functions, including its Bartholdi type, are 
	unified on a single scheme. 
	In Section 5, 
	we prepare some auxiliary facts relating the Schur complement of a matrix, 
	and verify the main theorem.

	Throughout this paper, we use the following notation. 
	The ring of integers is denoted by ${\mathbb Z}$, 
	and ${\mathbb Q}$ denotes the rational number field. 
	For a finite set $X$, the number of elements in $X$ is denoted by $|X|$, 
	and the family of subsets in $X$ by $2^X$. 
	The empty set is denoted by $\emptyset$. 
	For sets $X$ and $Y$, 
	$X\sqcup Y$ denotes the disjoint union of $X$ and $Y$. 
	The Kronecker delta is denoted by $\delta_{xy}$, 
	which gives one if $x=y$, zero if $x\neq y$. 
	For a proposition $P$, we also use the Kronecker product $\delta_P$, 
	which equals one if $P$ is true, zero otherwise. 
	For a square matrix $M$, the determinant of $M$ is denoted by $\det M$, 
	and the trace by ${\rm tr}\,M$. 
	For a ring $R$, $M_n(R)$ denotes the set of 
	square matrices of degree $n$ with entries in $R$. 
	The adjacency matrix and the degree matrix of a finite graph $\Gamma$ 
	are denoted by $A_{\Gamma}$ and $D_{\Gamma}$ respectively.

	\section{Preliminaries}
	\label{section : Preliminaries}
	
		\subsection{Digraphs}
		\label{subsection : Digraphs}
		
%
%
			\subsubsection{Digraphs and arcs}
			\label{subsubsection : Digraphs and arcs}
			
			A {\it digraph} $\Delta=(V, {\cal A})$ is a pair of a set $V$, 
			and a multi-set ${\cal A}$ consisting of ordered pair $(u, v)$ 
			of elements in $V$. 
			A digraph is {\it finite} if both $V$ and ${\cal A}$ are finite (multi-)sets. 
			An element of $V$ is called a {\it vertex}, 
			and ${\cal A}$ an {\it arc} respectively. 
			An arc $a=(u, v)$ is depicted by an arrow form $u$ to $v$. 
			For an arc $a=(u, v)$, 
			the vertex $u$ is called the {\it tail} of $a$, 
			and $v$ the {\it head} of $a$, 
			denoted by $\tail (a)$ and $\head(a)$ respectively. 
			Note that, since ${\cal A}$ is a multi-set, 
			it may occur that $\tail(a)=\tail(a')$ and $\head(a)=\head(a')$ for distinct $a, a'\in{\cal A}$. 
			An arc $l\in{\cal A}$ satisfying $\head(l)=\tail(l)$ is called a {\it loop} of $\Delta$. 
			For a loop $l$, 
			the vertex $n=\tail(l)(=\head(l))$ is called the {\it nest} of $l$. 
			We denote the set of loops by ${\cal L}$, 
			which is a subset of ${\cal A}$. 
			If we let 
			$$
				{\cal A}_{uv}=\{a\in{\cal A}\mid \tail(a)=u, \head(a)=v\}
			$$ 
			for $u, v\in V$, 
			then in general we have $|{\cal A}_{uv}|\geq 1$ if ${\cal A}_{uv}\neq\emptyset$. 
			If $u=v$, then ${\cal A}_{uu}$ consists of loops with nest $u$, 
			and we also allow $|{\cal A}_{uu}|\geq 1$ for each $u$. 
			A digraph is called {\it simple} 
			if ${\cal A}_{uv}\neq\emptyset$ implies $|{\cal A}_{uv}|=1$ 
			for  distinct $u, v\in V$, 
			and ${\cal A}_{uu}=\emptyset$ for every $u\in V$. 
			In general, 
			we do not suppose simplicity for digraphs in the present article. 
			Thus an element $a$ of ${\cal A}_{uv}$ is sometime called an {\it multi-arc}, 
			and the cardinality $|{\cal A}_{uv}|$ is called the {\it multiplicity} of $a$. 
			Similarly, an element of ${\cal A}_{uu}$ is called a {\it multi-loop} 
			with multiplicity $|{\cal A}_{uu}|$. 
			For $u, v\in V$, 
			let 
			\begin{equation}\label{equation : The arcs between u and v}
				{\cal A}(u, v)
				=
				{\cal A}_{uv}\cup {\cal A}_{vu}
			\end{equation}
			denotes the set 
			of arcs lying between the vertices $u$ and $v$. 
			If $u=v$, 
			then we have ${\cal A}(u, u)={\cal A}_{uu}$. 
			If $u$ and $v$ are distinct, 
			then the union 
			$
				{\cal A}(u, v)
				=
				{\cal A}_{uv}\sqcup {\cal A}_{vu}
			$ 
			is disjoint. 
			For $u, v\in V$, let 
			\begin{eqnarray*}
				&&
				{\cal A}_{u*}
				=
				\{a\in{\cal A}\mid \tail(a)=u\},\\
				&&
				{\cal A}_{*v}
				=
				\{a\in{\cal A}\mid \head(a)=v\}.
			\end{eqnarray*}

			\subsubsection{The inverses for an arc}
			\label{subsubsection : The inverses for an arc}

			Let $\Delta=(V, {\cal A})$ be a finite digraph, 
			and $a\in {\cal A}$ an arc. 
			If $a\in{\cal A}_{uv}$, 
			then any element $a'\in{\cal A}_{vu}$ is called an {\it inverse} of $a$, 
			that is, 
			${\cal A}_{vu}$ is the set of inverse arcs for $a\in{\cal A}_{uv}$. 
			Thus an inverse of an arc $a$ may not be uniquely determined for $a$. 
			We also note that any loop $l'\in{\cal A}_{uu}$ 
			is an inverse of $l\in{\cal A}_{uu}$, 
			containing $l$ itself. 
			In particular, each loop is {\it self-inverse}. 
			We denote ${\cal A}_{vu}$ by $S(a)$ if an arc $a$ belongs to ${\cal A}_{uv}$. 
			Note that $S(l)={\cal A}_{uu}$ for $l\in {\cal A}_{uu}$.

			This definition for inverse arcs is 
			a natural generalization for the usual one employed in the 
			theory of graph zeta functions in the case where 
			$\Delta$ is the symmetric digraph of a finite simple graph $\Gamma$. 
			To see this, we recall the usual definition of inverse arcs where $\Gamma$ is a finite graph, 
			not necessarily simple. 
			Let $\Gamma=(V, E)$ be a finite graph, 
			where $V$ and $E$ denote the vertex set and the edge set of $\Gamma$, respectively. 
			Namely, $E$ is a multi-set consisting of $2$-subsets $\{u, v\}$ of $V\times V$. 
			An element of $E$ is called a {\it multi-edge} in general. 
			In particular, if the multiplicity of an element $e\in E$ is one, 
			then $e$ is called a {\it simple edge}, or simply, an ${\it edge}$. 
			A multi-edge $\{u, u\}$ is called a {\it multi}-{\it loop} with {\it nest} $u$. 
			If $e=\{u, u\}$ is a simple edge, 
			then $e$ is called a {\it loop}. 
			Let $L$ denote the set of multi-loops of $\Gamma$. 
			If every element of $E$ is simple and $L=\emptyset$, 
			then $\Gamma$ is called a {\it simple graph}. 
			Let $\Gamma=(V, E)$ be a finite graph. 
			For each multi-edge $e=\{u, v\}\in E\setminus L$, 
			we assign two distinct arcs $a_e$ and $\overline{a_e}$ for $e$. 
			We set $\overline{\overline{a_e}}=a_e$, 
			and $a_e\neq a_{e'}$ for $e, e'\in E$ with $e\neq e'$. 
			It is assumed that the arc $a_e$ belongs to ${\cal A}_{uv}$ or ${\cal A}_{vu}$. 
			If $a_e\in {\cal A}_{uv}$, then we have $\overline{a_e}\in {\cal A}_{vu}$, 
			and vice versa. 
			If $e=\{u, u\}\in L$, 
			then we assign a single directed loop $l_e\in{\cal A}_{uu}$ for $e$. 
			If $e\neq e'$ in $L$, 
			then $l_e\neq l_{e'}$. 
			Let 
			$$
				{\cal A}=\{a_e, \overline{a_e}\mid e\in E\setminus L\}\sqcup\{l_e\mid e\in L\}. 
			$$
			The digraph $(V, {\cal A})$ is called the {\it symmetric digraph} of $\Gamma$, 
			denoted by $\Delta(\Gamma)$. 
			For the symmetric digraph $\Delta=\Delta(\Gamma)$ of a finite graph, 
			the usual definition of inverse arc is as follows 
			(see, e.g., \cite{bass92, bartholdi99, mizunosato04, sato07}). 
			Let $a\in{\cal A}$. 
			If $a=a_e$ for some $e\in E\setminus L$, 
			then $\overline{a_e}$ is the unique inverse arc for $a$, and vice versa. 
			If $a=l_e$ is a loop ($e\in L$), 
			then $a=l_e$ itself is the unique inverse arc for $a$. 
			

			\begin{lem}\label{lem : The arcs of the symmetric digraph for a finite simple graph} 
			Suppose that $\Delta=(V, {\cal A})$ is the symmetric digraph $\Delta(\Gamma)$ 
			of a finite simple graph $\Gamma$. Then we have 
				$
				{\rm i)}\; 
					{\cal A}_{uu}=\emptyset, \forall u\in V,\ 
				{\rm ii)}\; 
					{\cal A}_{uv}\neq\emptyset\Rightarrow |{\cal A}_{uv}|=1,\ 
				{\rm iii)}\; 
					|{\cal A}_{uv}|=1\Leftrightarrow|{\cal A}_{vu}|=1. 
			$
			\end{lem}
			
			By Lemma \ref{lem : The arcs of the symmetric digraph for a finite simple graph}, 
			one can easily see that 
			if ${\cal A}_{uv}\neq \emptyset$ then $|{\cal A}_{uv}|=|{\cal A}_{vu}|=1$ 
			for distinct $u, v\in V$. 
			This shows that those definitions of inverse arcs coincide with each other 
			on the symmetric digraph of a finite simple graph.

%
%

			\subsubsection{Closed paths}
			\label{subsubsection : Closed paths}
			
			Let $\Delta=(V, {\cal A})$ be a digraph. 
			Recall that $S(a)={\cal A}_{vu}$ for $a\in {\cal A}_{uv}$. 
			A sequence $c=(a_1, a_2, \dots , a_{m})$ of arcs is called a {\it path} 
			if it satisfies ${\frak h}(a_i)={\frak t}(a_{i+1})$ for each $i=1, 2, \dots , m-1$. 
			If $a_{i+1}\in S(a_i)$ fo some $i$, 
			then the pair $(a_i, a_{i+1})$ of arcs $a_i, a_{i+1}$ is called 
			a {\it backtrack} of $c$, 
			or sometimes called a {\it backtrack thorough} $\head(a_{i})$. 
			Then $m$ is called the {\it length} of $x$, denoted by $l(x)$. 
			A {\it closed} path is a path $c=(a_1, a_2, \dots , a_{m})$ with ${\frak h}(a_{m})={\frak t}(a_1)$. 
			A closed path $c=(a_1, a_2, \dots , a_{m})$ is called {\it reduced} 
			if it has no backtracks, i.e., $a_{i+1}\notin S(a_i)$ for any $i=1,2,\dots, m$, 
			where $a_{m+1}=a_1$. 

			Let $C_m=C_m(\Delta)$ denotes the set of closed paths of length $m$ in $\Delta$. 
			The set of reduced closed paths of length $m$ is denoted by $C^{\flat}_m=C^{\flat}_m(\Delta)$. 
			The set of closed paths in $\Delta$ is given by $C=C(\Delta):=\sqcup_{m\geq 1}C_m$. 
			The set of reduced closed paths is given by $C^{\flat}=C^{\flat}(\Delta):=\sqcup_{m\geq 1}C^{\flat}_m$. 
			For a positive integer $k$, 
			the concatenation of $k$ copies of a closed path $c=(a_1, a_2, \dots , a_{m})$ is also a closed path of $\Delta$. 
			This closed path is called the $k$-{\it th power} of $c$, 
			denoted by $c^k$. 
			The length of $c^k$ is $km$ if $l(c)=m$. 
			If there exists no shorter closed path $c'$ satisfying $c=c'^k$, 
			then $c$ is called {\it prime}. 
			We denote the set of prime closed paths of length $m$ by $P_m=P_m(\Delta)$. 
			The set of prime reduced closed paths of length $m$ is denoted by $P^{\flat}_m=P^{\flat}_m(\Delta)$. 
			The set of prime (resp. prime reduced) closed paths of $\Delta$ is given by 
			$
				P=P(\Delta):=\sqcup_{m\geq 1}P_m
			$ 
			(resp. 
			$
				P^{\flat}=P^{\flat}(\Delta):=\sqcup_{m\geq 1}P^{\flat}_m.
			$ 
			)

			\subsubsection{Cycles}
			\label{subsubsection : Cycles}

			The cyclic permutation $\sigma=(1,2,\dots ,m)$ acts on $C_m$ by 
			$$
				(a_1, a_2, \dots , a_{m}).\sigma
				=
				(a_{\sigma(1)}, a_{\sigma(2)}, \dots, a_{\sigma(m)}).
			$$ 
			Two closed paths $c=(a_i)$, $c'=(a_i')$ which belong to $C_m$ are called 
			{\it cyclically equivalent}, or simply {\it equivalent}, 
			if there exists an integer $k$ satisfying $c'=c.\sigma^k$. 
			The equivalence is denoted by $c\sim c'$. 
			The binary relation $\sim$ is indeed an equivalence relation. 
			An equivalence class $[c]=c\mbox{ mod }\sim$ is called an {\it cycle} of $\Delta$. 
			If $c\sim c'$ for $c, c'\in C$, 
			then we have $l(c)=l(c')$. 
			Hence one can define the {\it length} $l([c])$ of a cycle $[c]$ ($c\in C$) 
			by $l([c]):=l(c)$. 
			Let $[C]=C/\sim$ and $[C_m]=C_m\sim$ ($m\geq 1$). 			
			We have ${[C]}=\sqcup_{m\geq 1}{[C_m]}$. 
			If $c\in C$ is prime (resp. reduced), 
			then one can easily see that $c'\in C$ equivalent to $c$ is also prime (resp. reduced). 
			Hence we can say that a cycle $[c]$ is {\it prime} (resp. {\it reduced}) 
			if a representative $c$ is prime (resp. reduced). 
			A {\it prime reduced cycle} is a cycle $[c]$ with a representative $c$ which is prime and reduced. 
			The set of prime cycles, reduced cycles, prime reduced cycles are given by 
			$$
				{[P]}=P/\sim, 
				\quad
				{[C^{\flat}]}=C^{\flat}/\sim, 
				\quad
				{[P^{\flat}]}=P^{\flat}/\sim, 
			$$ 
			respectively. 
			These also have the decomposition 
			$
				{[P]}=\sqcup_{m\geq 1}[P_m], 
				{[C^{\flat}]}=\sqcup_{m\geq 1}[C^{\flat}_m], 
				{[P^{\flat}]}=\sqcup_{m\geq 1}[P^{\flat}_m], 
			$
			where $[P_m], [C^{\flat}_m], [P^{\flat}_m]$ denote 
			$P_m/\sim, C^{\flat}_m/\sim, P^{\flat}_m/\sim$ respectively. 
			Note that these are all well-defined.

			\subsubsection{The adjacency matrix and the backtrack matrix}
			\label{subsubsection : The adjacency matrix and the backtrack matrix}

			Let $\Delta=(V, {\cal A})$ be a finite digraph. 
			The matrix 
			$
				A_{\Delta}=(a_{uv})_{u, v\in V}
			$ 
			with entries $a_{uv}=|{\cal A}_{uv}|$ 
			is called the {\it adjacency matrix} of $\Delta$. 
			For each $u\in V$, 
			let $d_u$ denote the number of closed paths 
			$c=(a_1, a_2)$ satisfying i) $u=\tail(a_1)=\head(a_2)$, 
			ii)  $a_2\in S(a_1)$,  and iii) $a_1\notin{\cal L}$. 
			In other words, 
			$d_u$ is the number of backtracks though $u$. 
			The diagonal matrix $D_{\Delta}=(\delta_{uv}d_u)_{u, v\in V}$ 
			is called the {\it backtrack matrix} of $\Delta$.

			Note that the backtrack matrix $D_{\Delta}$ depends on the definition of inverse arcs. 
			Let $\Gamma=(V, E)$ be a finite graph, 
			and $\Delta=\Delta(\Gamma)=(V, {\cal A})$ the symmetric digraph of $\Gamma$. 
			For each $e=\{u, v\}\in E\setminus L$, we assign two arcs $a_e, \overline{a_e}$, 
			say $a_e\in{\cal A}_{uv}$. 
			Fix a vertex $u\in V$. 
			Note that if $\{u, v\}\notin E$ for $v\in V$, 
			then ${\cal A}_{uv}(={\cal A}_{vu})=\emptyset$. 
			If we employ the definition of inverse arcs as in \cite{bartholdi99, bass92} 
			(see \ref{subsubsection : The inverses for an arc}), 
			then there corresponds a unique backtrack $(\overline{a_e}, a_e)$ thorough $u$ 
			for each multi-edge $e$ of the form $\{u, v\}$ for some $v\in V$, $v\neq u$. 
			Thus we have $d_u=|{\cal A}_{u*}|$ in this case, 
			and one can readily see that this coincides with the degree 
			of the vertex $u$ in $\Gamma$. 
			Therefore, we have $D_{\Delta(\Gamma)}=D_{\Gamma}$ in this case. 
			On the other hand, 
			if we work on the definition of inverse arc introduced in \ref{subsubsection : The inverses for an arc}, 
			then any element of $A_{vu}$ is an inverse of $a_e$. 
			Hence, the number of backtracks through $u$ equals 
			the sum of $|{\cal A}_{uv}||{\cal A}_{vu}|=|{\cal A}_{uv}|^2$ for $v\in V$ with $\{u, v\}\in E$, 
			i.e., 
			$
				d_u
				=
				\sum_{v\in V}|{\cal A}_{uv}|^2.
			$
			If $\Gamma$ is a finite simple graph, 
			then it follows from Lemma \ref{lem : The arcs of the symmetric digraph for a finite simple graph} that 
			$d_u$ equals the degree of the vertex $u$ in $\Gamma$, 
			and we have $D_{\Delta(\Gamma)}=D_{\Gamma}$.

		\subsection{Words}\label{subsection : Words}

		\subsubsection{Definitions}
		\label{subsubsection : Definitions}
		Let ${\frak A}=\{\alpha_1, \alpha_2, \dots, \alpha_n\}$ be a finite alphabet, 
		and ${\frak A}^*$ the free monoid generated by $\frak A$. 
		An element of ${\frak A}^*$ is called a {\it word} on $\frak A$. 
		Let $w=a_1a_2\dots a_m$ be a word. 
		The integer $m$ is called the {\it length} of $w$, 
		denoted by $|w|$. 
		The multiplication on ${\frak A}^*$ is 
		defined by the concatenation of words. 
		Given two words $w, w'\in{\frak A}^*$, 
		the product of $w$ and $w'$ is denoted by $ww'$. 
		The $k$-th power of a word $w$ is denoted by $w^k$. 
		If a word $w$ can not be written by a power of a shorter word, 
		then $w$ is called a {\it prime word}. 
		Given a word $w=a_1a_2\cdots a_m\in{\frak A}^*$, 
		the {\it cyclic rearrangement class} ${\rm Re}\,(w)$ is the (multi-)set consisting of 
		the following $m$ words 
		$$
			a_1a_2\cdots a_{m-1}a_m, 
			a_2a_3\cdots a_ma_1, 
			\dots, 
			a_ma_1\cdots a_{m-2}a_{m-1}.
		$$
		If $w$ is a prime word, then its cyclic rearrangement class is a set.

		\subsubsection{Lyndon words}\label{subsubsection : Lyndon words}
		
		Let ${\frak A}=\{\alpha_1, \alpha_2, \dots, \alpha_n\}$ be a finite alphabet, 
		which is totally ordered by $\alpha_1<\alpha_2<\dots<\alpha_n$. 
		In this case, 
		the free monoid ${\frak A}^*$ is also totally ordered by 
		the lexicographical order induced by the total order $<$ on ${\frak A}$. 
		We denote the total order on ${\frak A}^*$ by the same symbol $<$. 
		If a word $w\in {\frak A}^*$ is the minimum element in its cyclic rearrangement class ${\rm Re}(w)$, 
		then $w$ is called a {\it Lyndon word} (see e.g., \cite{lothaire83}). 
		The set of Lyndon words on ${\frak A}$ is denoted by ${\rm Lyn}({\frak A})$. 
		For example, 
		if ${\frak A}=\{1<2<3\}$, 
		then $w=1212\notin {\rm Lyn}({\frak A})$, 
		since $w$ is not the minimum element in ${\rm Re}(w)=\{1212,2121,1212,2121\}$. 
		One can readily see from this example that 
		a Lyndon word is necessarily a prime word. 
		If $w=1213$, then $w$ is Lyndon. 
		The well-known Lyndon factorization theorem (c.f., \cite{lothaire83}) states that 
		the Lyndon words ${\rm Lyn}({\frak A})$ gives the primes of the free monoid ${\frak A}^*$.

		\subsubsection{The Foata-Zeilberger theorem}\label{subsubsection : The Foata-Zeilberger theorem}
		Let $R$ be a commutative ring, and 
		${\frak A}=\{\alpha_1<\alpha_2<\dots<\alpha_n\}$ a totally ordered finite alphabet. 
		Let ${\rm Mat}_{\frak A}(R)$ denote the set of $n\times n$ matrices $(m_{aa'})_{a, a'\in{\frak A}}$ with 
		$m_{aa'}\in R$ for each $a, a'\in{\frak A}$. 
		For $w=a_1a_2\cdots a_k\in{\frak A}^*$ and $M=(m_{aa'})_{a, a'\in{\frak A}}\in {\rm Mat}_{\frak A}(R)$, 
		let $\cir_M(w)$ denote the {\it circular product}
		$$
			m_{a_1a_2}m_{a_2a_3}\cdots m_{a_{k-1}a_k}m_{a_ka_1}
		$$
		of entries in $M$ along $w$. 
		Let $I$ denote the identity matrix of degree $n$ and $t$ an indeterminate. 
		The following proposition is called the {\it Foata-Zeilberger theorem}\cite{foatazeilberger99}. 
				
		\begin{prop}[Foata-Zeilberger]
			$	
			\det(I-M)
			=
			\prod_{l\in{\rm Lyn}(\frak A)}
			(1-\cir_M(l)).
			$
		\end{prop}

		It follows from the Foata-Zeilberger theorem that 
		the inverse of the reciprocal characteristic polynomial $1/\det(I-tM)$ 
		is written by 
		$$
			\prod_{l\in{\rm Lyn}(\frak A)}
			\frac{1}
			{1-\cir_M(l)t^{|l|}}.
		$$
		This identity can be viewed as the Euler product expression for $1/\det(I-tM)$ 
		since the set ${\rm Lyn}(\frak A)$ gives the primes in ${\frak A}^*$.

%
%

		\subsection{Dynamical systems}
		\label{subsection : Dynamical systems}
		
		\subsubsection{Prime period}
		\label{subsubsection : Prime period}
		
		A {\it dynamical system} is a pair $(\Xi, \lambda)$ of a set $\Xi$ and a bijection $\lambda : \Xi\rightarrow \Xi$. 
		For an positive integer $m$, 
		an element $x\in\Xi$ is called an $m$-{\it periodic} point of $(\Xi, \lambda)$ if 
		the condition $\lambda^m(x)=x$ holds. 
		The set of $m$-periodic points in $(\Xi, \lambda)$ is denoted by $X_m=X_m(\Xi)$, 
		and the set of all periodic points in $(\Xi, \lambda)$ by $X=X(\Xi)$, 
		We have $X=\cup_{m\geq 1}X_m$. 
		If $x\in X_m$, 
		then the positive integer $m$ is called a {\it period} of $x$. 
		Remark that any multiple of a period of $x$ is also a period of $x$. 
		Let ${\rm Per}(x)$ denote the set $\{m\mid x\in X_m\}$ of periods of $x\in\Xi$, 
		and we denote $\varpi(x)=\min {\rm Per}(x)$, 
		which is called the {\it prime period} of $x$. 
		It obviously follows that $x\in X_{\varpi(x)}$ for any $x\in X$. 
		If we regard $x$ as an element of $X_{\varpi(x)}$, 
		then we denote it by $\pi(x)$. 
		We call $\pi(x)$ the {\it prime section} of $x$. 

		If we denote by $Y_m=Y_m(\Xi)$ 
		the set of periodic points with prime period $m$, i.e., 
		$$
			Y_m=
			\{
				x\in X\mid \varpi(x)=m
			\},
		$$ 
		then we have a disjoint union $X=\sqcup_{m\geq 1}Y_m$. 
		A standard argument shows that $\varpi(x)$ divides any period of $x\in X$. 
		See, e.g., Lemma 8 in \cite{morita20} for a proof.

			\begin{lem}\label{lem:PrimeDividesPeriod}
			Let $x\in X$. We have $\varpi(x)|m$ for any $m\in{\rm Per}(x)$. 
			\end{lem}

		\subsubsection{Equivalence}
		\label{subsubsection : Equivalence}
		
		Two elements $x, y\in \Xi$ is called {\it equivalent} in $(\Xi, \lambda)$ 
		if there exists an integer $k$ satisfying $y=\lambda^k(x)$, 
		where $\lambda^{-1}$ denotes the inverse map of $\lambda$. 
		If $x$ and $y$ are equivalent, 
		we denote it by $x\equiv y$. 
		The set $\Xi/\equiv$ of equivalence classes is denoted by $[\Xi]$, 
		and an equivalence class $\xi\in[\Xi]$ with representative $x\in \Xi$ is denoted by $[x]$. 
		Given $x\in X$, the equivalence class $[x]$ is called the {\it orbit through} $x$. 
		The following lemma is verified in a straightforward way. 
		See e.g., Lemma 4 in \cite{morita20}.

			\begin{lem}
			\label{lem : Periods Equal}
			If $x\equiv y$ in $(\Xi, \lambda)$, 
			then ${\rm Per}(x)={\rm Per}(y)$. 
			Hence we have $\varpi(x)=\varpi(y)$. 
			\end{lem}
%

		Lemma \ref{lem : Periods Equal} shows that the equivalence relation $\equiv$ is an equivalence relation on $X$, 
		and also on each $X_m$. 
		The equivalence classes $X/\equiv$ (resp. $X_m/\equiv$) is denoted by 
		${[X]}$ (resp. ${[X]}_m$). 
		We have ${[X]}=\cup_{m\geq 1}{[X]}_m$. 
		Note that, by Lemma \ref{lem : Periods Equal}, 
		the binary relation $\equiv$ is also an equivalence relation on $Y_m$ for each $m$. 
		We denote $Y_m/\equiv$ by ${[Y_m]}$. 
		It also suggests that 
		we can define the {\it prime period} $\varpi([x])$ of an orbit $[x]\in{[X]}$ by 
		$\varpi([x])=\varpi(x)$. 
			

		\subsubsection{Ruelle zeta functions}
		\label{subsubsection : Ruelle zeta functions}
		
		Let $(\Xi, \lambda)$ be a dynamical system, 
		and $X=\cup_{m\geq 1}X_m$ the set of periodic points in $(\Xi, \lambda)$, 
		where $X_m$ denotes the set of $m$-periodic points as in \ref{subsubsection : Prime period}. 
		If each $X_m$ is a finite set, 
		then the dynamical system $(\Xi, \lambda)$ is called {\it quasi-finite}. 
		Let $R$ be a commutative ${\mathbb Q}$-algebra, 
		and let $\chi_m : X_m\rightarrow R$ be a map. 
		The symbol $\chi$ denotes the multi-valuated map $X\rightarrow R$ 
		which sends $x\in X_m$ to $\chi_m(x)$. 
		The triple $(\Xi, \lambda, \chi)$ is called a {\it weighted dynamical system} with a weight $\chi$ (c.f., \cite{hattorimorita16}). 
		If a weighted dynamical system $(\Xi, \lambda, \chi)$ is quasi-finite, 
		then the following sum 
		$
			N_m(\chi)
			=
			\sum_{x\in X_m}\chi(x)
		$ 
		is well-defined for each $m\geq 1$.

		\begin{df}[c.f., \cite{ruelle94}]
		{\em 
		Let $t$ be an indeterminate, 
		and $(\Xi, \lambda, \chi)$ a quasi-finite weighted dynamical system. 
		The {\it Ruelle zeta function} for $(\Xi, \lambda, \chi)$ is 
		the formal power series 
		\begin{equation}\label{equation : The Ruelle zeta function}
			\exp
				\left(
					\sum_{m\geq 1}
					\frac{N_m(\chi)}{m}
					t^m
				\right),
		\end{equation}
		which we denote by $Z_{\Xi}(t; \chi)$. 
		The Ruelle zeta function $Z_{\Xi}(t; \chi)$ is also simply called 
		the {\it zeta function} for $(\Xi, \lambda, \chi)$.
		}
		\end{df}

		\section{Zeta functions for digraphs}
		\label{section : Zeta functions for digraphs}
		
		In this section, 
		we briefly review the theory of combinatorial zeta functions, 
		which is the fundamental framework for considering the Hashimoto expression. 
		See \cite{morita20} for precise information on the development in this section. 
		
		\subsection{Dynamical systems on digraphs}
		\label{subsection : Dynamical systems on digraphs}
		
		Let $\Delta=(V, {\cal A})$ be a finite digraph. 
		We denote by 
		$
			{{\cal A}}^{\mathbb Z}
		$ 
		the set 
		$
			\{(a_i)_{i\in{\mathbb Z}}\mid a_i\in{A}, \forall i\in{\mathbb Z}
			\}
		$
		of two-sided infinite sequences with entries in ${\cal A}$. 
		Consider the left shift operator 
		$$
			\varphi : {{\cal A}}^{\mathbb Z}\rightarrow {{\cal A}}^{\mathbb Z} : 
			(a_i)_i\mapsto(a_{i+1})_i
		$$ 
		on ${{\cal A}}^{\mathbb Z}$. 
		Let $x=(x_i)_{i\in{\mathbb Z}}\in {\cal {\cal A}}^{\mathbb Z}$ and let $m$ be a positive integer. 
		Any finite consecutive segment $(x_i, x_{i+1}, \dots, x_{i+m-1})$ of $x$ is called 
		an $m$-{\it section}. 
		In particular, 
		the $m$-section $(x_0, x_1, \dots, x_{m-1})$ is called 
		the {\it principal} $m$-{\it section} of $x$, denoted by ${\rm ps}_m(x)$. 
		A subset $\Xi\subset {{\cal A}}^{\mathbb Z}$ is called $\varphi$-{\it stable} 
		if it satisfies the condition $\varphi(\Xi)\subset\Xi$.

		Suppose that $\Xi\subset {{\cal A}}^{\mathbb Z}$ is $\varphi$-{stable}. 
		We denote the restriction $\varphi|_{\Xi}$ by $\lambda$, 
		and consider the dynamical system $(\Xi, \lambda)$. 
		Recall that $S(a)$ denotes ${\cal A}_{vu}$ for $a\in{\cal A}_{uv}$. 
		The examples that we mainly bare in mind are the set $\Pi_{\Delta}$ of two-sided infinite paths of $\Delta$ 
		and the set $\Pi^{\flat}_{\Delta}$ of two-sided infinite reduced paths of $\Delta$, 
		i.e., 
		$$
			\begin{array}{l}
			\Pi_{\Delta}=
				\{(a_i)_{i\in{\mathbb Z}}\in {\cal {\cal A}}^{\mathbb Z}
				\mid \head(a_i)=\tail(a_{i+1}),\ \forall i\},\\
			\Pi^{\flat}_{\Delta}=
				\{(a_i)_{i\in{\mathbb Z}}\in {\cal {\cal A}}^{\mathbb Z}
				\mid \head(a_i)=\tail(a_{i+1}),\ a_{i+1}\notin S(a_i),\ \forall i\}. 
			\end{array}
		$$
		Consider the case where $\Xi=\Pi_{\Delta}$. 
		Note that the dynamical system $(\Xi, \lambda)$ is quasi-finite 
		since the cardinality $|X_m|$ of $m$-periodic points does not exceed $|A|^m$. 
		Note also that an $m$-section of $x\in X_m$ corresponds to an element of $C_m=C_m(\Delta)$. 
		Since $x\in X_m$ is completely determined by an $m$-section, 
		in particular by the principal $m$-section ${\rm ps}_m(x)$, 
		the map 
		$$
			{\rm ps}_m : X_m\rightarrow C_m : x\mapsto {\rm ps}_m(x)
		$$ 
		is bijective. 
		For $c=(a_0, a_1, \dots, a_{m-1})\in C_m$, 
		we denote by $c^{\natural}$ the element $x=(x_i)_{i\in{\mathbb Z}}\in X_m$ defined by 
		$
			x_i=a_k
		$ for $i\in{\mathbb Z}$ with the indices $i$ congruent with $k$ modulo $m$. 
		It is clear that the map 
		$$
			\natural : C_m\rightarrow X_m : c\mapsto c^{\natural}
		$$ 
		gives the inverse for the bijection ${\rm ps}_m$. 
		Suppose that $x\in\Xi$ belongs to $Y_m$, the set of elements with prime period $m$. 
		Since $m=\min {\rm Per}(x)$, 
		the image ${\rm ps}_m(x)$ belongs to $P_m=P_m(\Delta)$. 
		Thus the restriction ${\rm ps}_m|_{Y_m}$ gives a bijective correspondence $Y_m\rightarrow P_m$. 
		In the case where $\Xi=\Pi_{\Delta}^{\flat}$, 
		we denote the set of $m$-periodic points by $X_m^{\flat}=X_m^{\flat}(\Delta)$. 
		The set of points in $X_m^{\flat}$ with prime period $m$ is denoted by $Y_m^{\flat}=Y_m^{\flat}(\Delta)$. 
		In the same manner, 
		one can see that 
		$X_m^{\flat}$ (resp. $Y_m^{\flat}$) is mapped bijectively onto $C_m^{\flat}$ (resp. $P_m^{\flat}$) 
		by the map ${\rm ps}_m$.

		Let $\Xi\subset {\cal {\cal A}}^{\mathbb Z}$ be $\varphi$-stable, 
		and consider the dynamical system $(\Xi, \lambda)$ where $\lambda=\varphi|_{\Xi}$. 
		Let $\equiv$ denote the equivalence relation on $(\Xi, \lambda)$ introduced in 
		\ref{subsection : Dynamical systems}. 
		Note that, if $x\equiv x'$ for $x, x'\in\Xi_m$, 
		then it follows that ${\rm ps}_m(x)\sim{\rm ps}_m(x')$. 
		This shows that the map ${\rm ps}_m : X_m\rightarrow C_m$ induces 
		a  bijection 
		$$
			[{\rm ps}_m] : [X_m]\rightarrow [C_m] : [x]\mapsto [{\rm ps}_m(x)]
		$$
		for each $m\geq 1$. 
		It also follows that the restrictions $[{\rm ps}_m] |_{[Y_m]} : [Y_m]\rightarrow[P_m]$, 
		$[{\rm ps}_m]|_{[X_m^{\flat}]} : [X_m^{\flat}]\rightarrow [C_m^{\flat}]$, 
		$[{\rm ps}_m]|_{[Y_m^{\flat}]} : [Y_m^{\flat}]\rightarrow [P_m^{\flat}]$ are bijective.

		Let $[X]:=X/\equiv$ and $[X^{\flat}]:=X^{\flat}/\equiv$. 
		Since $\equiv$ is an equivalence relation on each $X_m$ and $X_m^{\flat}$, 
		one has $[X]=\cup_m[X_m]$ and $[X^{\flat}]=\cup_m[X_m^{\flat}]$. 
		Let $[Y]:=Y/\equiv$ and $[Y^{\flat}]:=Y^{\flat}/\equiv$. 
		Since $\equiv$ is an equivalence relation on each $Y_m$ and $Y_m^{\flat}$, 
		one has $[Y]=\sqcup_m[Y_m]$ and $[Y^{\flat}]=\sqcup_m[Y_m^{\flat}]$. 
		It also follows that $[Y]=\pi([X])$ and $[Y^{\flat}]=\pi([X^{\flat}])$.

		\subsection{The path condition}
		\label{section : The path condition}
		
			Let $\Delta=(V, {\cal A})$ be a finite digraph, 
			$S(a)={\cal A}_{vu}$ for $a\in{\cal A}_{uv}$, 
			$\Xi$ a $\varphi$-stable subset of ${\cal A}^{\mathbb Z}$, 
			and $\theta : {\cal A}\times{\cal A}\rightarrow R$ a map. 
			Recall that 
			$w^{\natural}$ is the element $(x_i)_{i\in{\mathbb Z}}$ of ${\cal A}^{\mathbb Z}$ 
			given by $x_i=a_j$ if $i$ is congruent to $j$ modulo $k$. 
			We give a total order ${\cal A}=\{\alpha_1<\alpha_2<\cdots<\alpha_n\}$ on ${\cal A}$, 
			and consider the set ${\rm Lyn}(\cal A)$ of Lyndon words on the alphabet ${\cal A}$. 
			For a word $w=a_0a_1\cdots a_{k-1}\in{\cal A}^*$, 
			we denote the circular product 
			$$
				\theta(a_0, a_1)\theta(a_1, a_2)\cdots \theta(a_{k-2}, a_{k-1})\theta(a_{k-1}, a_1)
			$$ 
			by $\cir_{\theta}(w)$. 
			The following {\lq\lq path condition\rq\rq} was introduced in 3.2.1 in \cite{morita20}.

			\begin{df}[The path condition]
			\label{def : The path condition}
			{\em
			We say that a $\varphi$-stable subset $\Xi$ of ${\cal A}^*$ satisfies the {\it path condition} with respect to $\theta$ 
			iff one has 
			$
				\cir_{\theta}(l)\neq 0
				\Rightarrow 
				l^{\natural}\in\Xi
			$
			for each $l\in{\rm Lyn}({\cal A})$. 
			}
			\end{df}

			In this case, we also say that $(\Xi, \theta)$ satisfies the path condition. 
			We will see some examples. 
			Let 
			$
				\theta^{\rm BL}
			$ and 
			$\theta^{{\rm BL}, \flat}
			$
			be maps ${\cal A}\times {\cal A}\rightarrow R$ given by 
			$\theta^{\rm BL}(a, a')=\delta_{\head(a)\tail(a')},$ 
			and 
			$\theta^{{\rm BL}, \flat}(a, a')=\delta_{\head(a)\tail(a')}-\delta_{a'\in S(a)}.$ 
			One can see that $(\Pi_{\Delta}, \theta^{\rm BL})$ satisfies the path condition, 
			since, for $l=a_0\cdots a_{m-1}\in{\rm Lyn}(\cal A)$, 
			$\cir_{\theta^{\rm BL}}(l)\neq 0$ implies 
			$\head(a_i)=\tail(a_{i+1})$ for any $i=0, \dots, m-1$ 
			where $a_m:=a_0$, 
			and hence it follows that $l^{\natural}\in\Pi_{\Delta}$. 
			It can also be shown that $(\Pi_{\Delta}^{\flat}, \theta^{{\rm BL}, \flat})$ 
			satisfies the path condition. 
			On the other hand, one can see that $(\Pi_{\Delta}^{\flat}, \theta^{\rm BL})$ 
			does not satisfy the path condition. 
			In particular, 
			if we consider the case where $\Xi=\Pi_{\Delta}$ or $\Pi_{\Delta}^{\flat}$ 
			as in these examples, 
			then it can be seen that 
			the path condition is implied by the following simple condition. 
			See 4.1 in \cite{morita20} for precise information. 
			
			\begin{df}
			\label{def : The adjacency condition}
			{\em 
			A map $\theta : {\cal A}\times{\cal A}\rightarrow R$ is said to satisfy 
			the {\it adjacency condition} iff 
			$\theta(a, a')\neq 0$ implies $\head(a)=\tail(a')$ for $a, a'\in{\cal A}$. 
			A map $\theta : {\cal A}\times{\cal A}\rightarrow R$ is said to satisfy 
			the {\it reduced adjacency condition} iff 
			$\theta(a, a')\neq 0$ implies $\head(a)=\tail(a')$ and $a'\notin S(a)$ for $a, a'\in{\cal A}$. 
			}
			\end{df}
			
			
			\begin{lem}
			\label{lem : adjacency conditions}
			If a map $\theta : {\cal A}\times{\cal A}\rightarrow R$ satisfies the adjacency condition, 
			then $(\Pi_{\Delta}, \theta)$ satisfies the path condition. 
			If $\theta$ satisfies the reduced adjacency condition, 
			then $(\Pi_{\Delta}^{\flat}, \theta)$ satisfies the path condition. 
			\end{lem}

		\subsection{Zeta functions for finite digraphs}
		\label{subsection : Zeta functions for finite digraphs}
		
		Let $\Delta=(V, {\cal A})$ be a finite digraph with $n$ arcs , 
		$S(a)={\cal A}_{vu}$ for $a\in{\cal A}_{uv}$, 
		$R$ a commutative ${\mathbb Q}$-algebra, 
		and 
		$\theta : {\cal A}\times{\cal A}\rightarrow R$ a map. 
		Let $\varphi$ be the left shift operator on ${\cal {\cal A}}^{\mathbb Z}$, 
		$\Xi$ a $\varphi$-stable subset of ${\cal {\cal A}}^{\mathbb Z}$, 
		and $\lambda$ the restriction $\varphi|_{\Xi}$. 	

			\subsubsection{Definition and the exponential expression}
			\label{subsubsection : Definition and the exponential expression}
			
			Let $\Delta=(V, {\cal A})$ be a finite digraph, 
			and $\Xi$ a $\varphi$-stable subsets of ${\cal A}^{\mathbb Z}$. 
			Given an $m$-periodic point $x=(a_i)_i\in X_m=X_m(\Xi)$, 
			$\chi(x)$ stands for the following product 
			$$
				\theta(a_1, a_2)\cdots \theta(a_{m-1}, a_m)\theta(a_m, a_1), 
			$$ 
			which is called the {\it circular product} of $\theta$ along $x\in X_m$, 
			and we denote it by $\cir_\theta(x)$. 
			Note that $\chi(x)$ depends on the choice of a period of $x$, 
			hence $\cir_\theta$ is not a map in general. 
			This multi-valuated map
			$
				\cir_{\theta} : X\rightarrow R
			$ 
			is called the {\it circular weight} induced by $\theta$. 
			Let $\chi=\cir_{\theta}$, 
			$
				N_m(\chi)
				=
				\sum_{x\in X_m}\chi(x),
			$ 
			and $t$ an indeterminate. 
			
			\begin{df}[Graph zeta functions]
			\label{def : Graph zeta functions}
			{\em 
			A {\it graph zeta function} is 
			the Ruelle zeta function 
			\begin{equation}
			\label{equation : The exponential expression}
				\exp
					\left(
						\sum_{m\geq 1}
							\frac{N_m(\chi)}{m}t^m
					\right)
			\end{equation}
			with the circular weight $\chi=\cir_{\theta}$, 
			which we denoted by 
			$
				Z_{\Xi}(t; \theta). 
			$ 
			}
			\end{df}
			
			The fundamental literature on graph zetas is Terras \cite{terras11}. 
			This defining identity (\ref{equation : The exponential expression}) is called the {\it exponential expression}. 
			In particular, 
			if $\Xi=\Pi_{\Delta}, \Pi_{\Delta}^{\flat}$, 
			then we denote $Z_{\Xi}(t; \theta)$ by 
			$Z_{\Delta}(t; \theta)$, $Z_{\Delta}^{\flat}(t; \theta)$ respectively. 
			In particular, if $\Delta$ is the symmetric digraph of a finite graph $\Gamma$, 
			then $Z_{\Delta}(t; \theta)$, $Z_{\Delta}^{\flat}(t; \theta)$ are denoted by 
			$Z_{\Gamma}(t; \theta)$, $Z_{\Gamma}^{\flat}(t; \theta)$ respectively. 
			In the following of this paper, 
			$N_m(\chi)$ is denoted by $N_m(\theta)$ if $\chi=\cir_{\theta}$.

			\subsubsection{The Euler expression}
			\label{subsubsection : The Euler expression}

			It is shown in \cite{morita20} that 
			the exponential expression $Z_{\Xi}(t; \theta)$ of a graph zeta can always be 
			reformulated into the \lq\lq Euler expression\rq\rq. 
			Recall that $\varpi(x)$ is the prime period and 
			$\pi(x)$ the prime section for a periodic element $x\in X$ (see \ref{subsubsection : Prime period}). 
			We consider the following formal power series 
			\begin{equation}\label{equation : The Euler expression}
				\prod_{[x]\in[X]}
				\frac{1}{1-\cir_{\theta}(\pi(x))t^{\varpi(x)}},
			\end{equation}
			which is denoted by $E_{\Xi}(t; \theta)$. 
			The following proposition is verified in \cite{morita20}. 
			
			\begin{prop}[\cite{morita20}, Theorem 18]
			\label{prop : The Euler expression for graph zetas}
				$Z_{\Xi}(t; \theta)=E_{\Xi}(t; \theta).$
			\end{prop}
			This identity is called the {\it Euler expression} of the graph zeta function $Z_{\Xi}(t; \theta)$. 
			Note that, as one can see in \cite{morita20}, 
			the Euler expression does not need the path condition for $(\Xi, \theta)$.

			The zeta function $Z_{\Xi}(t; \theta)$ is 
			an example of \lq\lq zeta function associated with a family of finite set\rq\rq 
			introduced in \cite{morita20}. 
			By the notation in \cite{morita20}, 
			a graph zeta $Z_{\Xi}(t; \theta)$ is denoted by $Z_{\cal F}(t; \chi)$, 
			where 
			${\cal F}=\{X_m=X_m(\Xi)|m\geq 1\}$, $\chi=\cir_{\theta}$. 
			Since $({\cal F}, \chi, \equiv)$ satisfies the {\it Euler condition} (see 3.1.2, \cite{morita20}), 
			the zeta $Z_{\cal F}(t; \chi)$ has the Euler expression $E_{\cal F}(t; \chi)$ 
			(see 3.1.2 of \cite{morita20} for $E_{\cal F}(t; \chi)$), 
			and one can confirm that this coincides with $E_{\Xi}(t; \theta)$ above. 

			In the case where $\Xi=\Pi_{\Delta}$ (resp. $\Pi_{\Delta}^{\flat}$), 
			we denote $E_{\Xi}(t; \chi)$ by 
			$E_{\Delta}(t; \theta)$ (resp. $E^{\flat}_{\Delta}(t; \theta)$) 
			for $\chi=\cir_{\theta}$. 
			For the symmetric digraph $\Delta=\Delta(\Gamma)$ for a finite graph $\Gamma$, 
			$E_{\Delta}(t; \theta)$ (resp. $E^{\flat}_{\Delta}(t; \theta)$) is denoted by 
			$E_{\Gamma}(t; \theta)$ (resp. $E^{\flat}_{\Gamma}(t; \theta)$). 
			See 3.2.2 of \cite{morita20} for further information.

			\subsubsection{The Hashimoto expression}
			\label{subsubsection : The Hashimoto expression}
			
			Suppose that $(\Xi, \theta)$ satisfies the path condition. 
%
			Fixing a total order on ${\cal A}$, 
			we consider the matrix $M=(\theta(a, a'))_{a, a'\in{\cal A}}$ 
			determined by the map $\theta$. 
			Since $(\Xi, \theta)$ satisfies the path condition, 
			it follows that 
			$
				\cir_{M}(l)\neq 0
			$ 
			implies 			$
				l^{\natural}\in \Xi
			$ for $l\in{\rm Lyn}({\cal A})$ 
			(c.f., \ref{subsubsection : The Foata-Zeilberger theorem} for $\cir_M(l)$). 
			This suggest that any closed cycle $[c]=[(a_1, \dots, a_m)]$ of $\Delta$ 
			satisfying $\cir_{\theta}([c])\neq 0$ indeed appears in $\Xi$. 
			We denote the matrix $M$ by $M_{\Xi}(\theta)$, 
			and call it the {\it edge matrix} for $(\Xi, \theta)$. 
			Now we denote the following formal power series
			$$
				\frac{1}{\det(I-tM_{\Xi}(\theta))}. 
			$$
			by $H_{\Xi}(t; \theta)$, 
			where $I$ stands for the identity matrix of degree $n=|{\cal A}|$. 

			\begin{prop}[\cite{morita20}, Theorem 19]
			\label{prop : The Hashimoto expression for graph zetas}
			If $(\Xi, \theta)$ satisfies the path condition, 
			then we have 
			$Z_{\Xi}(t; \theta)=H_{\Xi}(t; \theta).$
			\end{prop}
			This identity in Proposition \ref{prop : The Hashimoto expression for graph zetas} 
			is called the {\it Hashimoto expression} of a graph zeta function $Z_{\Xi}(t; \theta)$. 
			Contrary to the Euler expression, 
			the Hashimoto expression does need the path condition for $(\Xi, \theta)$ (c.f., \cite{morita20}). 
			In the case where $\Xi=\Pi_{\Delta}$ (resp. $\Pi^{\flat}_{\Delta}$), 
			$H_{\Xi}(t; \theta)$ is denoted by $H_{\Delta}(t; \theta)$ (resp. $H^{\flat}_{\Delta}(t; \theta)$). 
			In particular, 
			if $\Delta=\Delta(\Gamma)$ for a finite graph $\Gamma$, 
			then $H_{\Delta}(t; \theta)$ (resp. $H^{\flat}_{\Delta}(t; \theta)$) is denoted by 
			$H_{\Gamma}(t; \theta)$(resp. $H^{\flat}_{\Gamma}(t; \theta)$. 
			See 3.2.3 in \cite{morita20} for precise information.

	\section{The generalized weighted zeta function}
	\label{section : The generalized weighted zeta function}
		
	Let $\Delta=(V, {\cal A})$ be a finite digraph, 
	$S(a)={\cal A}_{vu}$ for $a\in {\cal A}_{uv}$, and 
	$R$ a commutative ${\mathbb Q}$-algebra with unity $1$. 
	We consider two functions $\tau$, $\upsilon$ from ${\cal A}$ to $R$, 
	and 
	let 
	$
		\theta^{\rm G} : {\cal A}\times{\cal A}\rightarrow R
	$ 
	be the map defined by 
	\begin{equation}\label{equation : The generalized weighted map}
		\theta^{\rm G}(a, a')
		=
		\tau(a')\delta_{\head(a)\tail(a')}
		-
		\upsilon(a')\delta_{a'\in S(a)}.
	\end{equation}

		\subsection{Definition and the three expressions}
		\label{subsection :　Definition and the three expressions}
		
		\begin{df}[The generalized weighted zeta]
		\label{def : The generalized weighted zeta function}
		{\em 
		Let $\Xi$ be a $\varphi$-stable subset of ${\cal A}^{\mathbb Z}$. 
		The formal power series $Z_{\Xi}(t; \theta^{\rm G})$ is called 
		the {\it generalized weighted zeta function}.
		}
		\end{df}
		
		The generalized weighted zeta was introduced in 4.4 of \cite{morita20}. 
		It follows from Proposition \ref{prop : The Euler expression for graph zetas} that 
		$
			Z_{\Xi}(t; \theta^{\rm G})
			=
			E_{\Xi}(t; \theta^{\rm G}), 
		$ 
		and the path condition for $(\Xi, \theta^{\rm G})$ implies 
		$
			E_{\Xi}(t; \theta^{\rm G})
			=
			H_{\Xi}(t; \theta^{\rm G})
		$ 
		by Proposition \ref{prop : The Hashimoto expression for graph zetas}. 
		Hence we have the three expressions for the generalized weighted zeta function $Z_{\Xi}(t; \theta^{\rm G})$. 
		
		\begin{prop}
		\label{prop : The three expressions for the generalized weighted zeta function} 
			If $(\Xi, \theta^{\rm G})$ satisfies the path condition, 
			then we have 
			$
				Z_{\Xi}(t; \theta^{\rm G})
				=
				E_{\Xi}(t; \theta^{\rm G})
				=
				H_{\Xi}(t; \theta^{\rm G}).
			$ 
		\end{prop}
		
		In particular for $\Xi=\Pi_{\Delta}$, $\Pi_{\Delta}^{\flat}$, 
		the generalized weighted zeta $Z_{\Xi}(t; \theta^{\rm G})$ is denoted by 
		$Z_{\Delta}(t; \theta^{\rm G})$, $Z^{\flat}_{\Delta}(t; \theta^{\rm G})$, respectively. 
		The pairs $(\Pi_{\Delta}, \theta^{\rm G})$ and $(\Pi^{\flat}_{\Delta}, \theta^{\rm G})$ 
		may, indeed, satisfy the path condition, 
		that is, 
		we have the following lemme, 
		which is verified in 4.4 of \cite{morita20}. 
		
		\begin{lem}
		\label{lem : The adjacency condition for G}
		The map $\theta^{\rm G}$ satisfies the adjacency condition. 
		If $\tau=\upsilon$, then $\theta^{\rm G}$ satisfies the reduced adjacency condition.
		\end{lem}

		\begin{cor}\label{cor : The three expression for the generalized weighted zeta}
		For a finite digraph $\Delta$, we have 
		$
			Z_{\Delta}(t; \theta^{\rm G})
			=
			E_{\Delta}(t; \theta^{\rm G})
			=
			H_{\Delta}(t; \theta^{\rm G}). 
		$ 
		In the case where $\tau=\upsilon$, we also have 
		$
			Z^{\flat}_{\Delta}(t; \theta^{\rm G})
			=
			E^{\flat}_{\Delta}(t; \theta^{\rm G})
			=
			H^{\flat}_{\Delta}(t; \theta^{\rm G}).
		$ 
		\end{cor}

		\subsection{Unification of graph zeta functions}
		\label{subsection : Unification of graph zeta functions}

		The notation $\Delta=(V, {\cal A})$, $R$, $\tau$, $\upsilon$ and $\theta^{\rm G}$ are inherited. 
		Let $\Xi$ be a $\varphi$-stable subset of ${\cal A}^{\mathbb Z}$, 
		where $\varphi$ is the left shift operator on ${\cal A}^{\mathbb Z}$. 
		The typical examples of $\Xi$ are $\Pi_{\Delta}$ and $\Pi_{\Delta}^{\flat}$. 
		Varying $\tau$ and $\upsilon$, 
		the generalized weighted zeta $Z_{\Xi}(t; \theta^{\rm G})$ 
		degenerates to various graph zetas which have been studied in previous research 
		\cite{bartholdi99, bowenlanford70, ihara66, mizunosato04, sato07, starkterras96} etc. 
		
			\subsubsection{The Ihara zeta function}
			\label{subsubsection : The Ihara zeta function}
			In the case where $\tau=\upsilon=1$, 
			the generalized weighted zeta $Z^{\flat}_{\Delta}(t; \theta^{\rm G})$ is called the {\it Ihara zeta function} of $\Delta$. 
			In this case,
			we denote the map $\theta^{\rm G}$ by $\theta^{\rm I}$. 
			By Lemma \ref{lem : The adjacency condition for G}, 
			the map $\theta^{\rm I}$ satisfies the reduced adjacency condition, 
			and hence the adjacency condition. 
			Thus both $(\Pi^{\flat}_{\Delta}, \theta^{\rm I})$ and $(\Pi_{\Delta}, \theta^{\rm I})$ 
			satisfy the path condition by Lemma \ref{lem : adjacency conditions}. 
			Hence we have the identities 
			$
				Z^{\flat}_{\Delta}(t; \theta^{\rm I})
				=
				E^{\flat}_{\Delta}(t; \theta^{\rm I})
				=
				H^{\flat}_{\Delta}(t; \theta^{\rm I})
			$
			and 
			$
				Z_{\Delta}(t; \theta^{\rm I})
				=
				E_{\Delta}(t; \theta^{\rm I})
				=
				H_{\Delta}(t; \theta^{\rm I}).
			$
			Since the circular weight $\cir_{\theta^{\rm I}}$ instinctively excludes the 
			non-reduced paths, 
			that is, 
			$\cir_{\theta^{\rm I}}(x)\neq 0$ impies $x\in\Pi_{\Delta}^{\flat}$, 
			it follows that 
			$
				N_m(\theta^{\rm I})
				=
				|X_m^{\flat}|
			$
			for both $\Pi_{\Delta}$ and $\Pi^{\flat}_{\Delta}$, 
			and we have 
			$
				Z_{\Delta}(t; \theta^{\rm I})
				=
				Z^{\flat}_{\Delta}(t; \theta^{\rm I}). 
			$
			Thus these six formal power series all equals.

			The Ihara zeta was originally defined for 
			the case where $\Delta=\Delta(\Gamma)$ is the symmetric digraph of a finite graph $\Gamma=(V, E)$ 
			\cite{bass92, ihara66} (see also \cite{hashimoto90, serre80, sunada88}). 
			The {\it Bass-Ihara theorem} \cite{bass92, ihara66} states that 
			the reciprocal of the Ihara zeta for $\Delta=\Delta(\Gamma)$ equals the following polynomial 
			$$
				(1-t^2)^{|E|-|V|}
				\det(I-tA_{\Gamma}+t^2(D_{\Gamma}-I)).
			$$
			This expression is called the {\it Ihara expression}  for $Z_{\Gamma}(t; \theta^{\rm I})$, 
			denoted by $I_{\Gamma}(t; \theta^{\rm I})$. 
			Thus the Ihara zeta has four expressions 
			$Z_{\Gamma}(t; \theta^{\rm I})$, $E_{\Gamma}(t; \theta^{\rm I})$,$H_{\Gamma}(t; \theta^{\rm I})$, 
			$I_{\Gamma}(t; \theta^{\rm I})$, 
			and this is the case for the other graph zetas which have been appeared in the previous studies,  
			for instance \cite{bartholdi99, mizunosato04, sato07}, as we will see in the sequel 
			(except the Bowen-Lanford zeta function in \ref{subsubsection : Bowen-Lanford zeta}). 
			However, our fundamental point of view in the present article is that 
			the graph zeta functions should be defined for digraphs. 
			The reason is that the Ihara expression, one of the main interest of the preceding studies 
			including the present article, 
			can be constructed for any finite digraph as we will see in our main theorem. 
			Therefore, although the original definitions of those traditional graph zetas were for finite graphs, 
			the definitions of them in this paper are described for finite digraphs.

			\subsubsection{The Bowen-Lanford zeta function}
			\label{subsubsection : Bowen-Lanford zeta}
			In the case where $\upsilon=0$, 
			the generalized weighted zeta $Z_{\Delta}(t; \theta^{\rm G})$ is called 
			the {\it weighted Bowen-Lanford zeta function} \cite{bowenlanford70} for $\Delta$ 
			(see also \cite{morita20}). 
			We denote the map $\theta^{\rm G}$ by $\theta^{{\rm BL}}_{\tau}$. 
			Note that the original definition in \cite{bowenlanford70} is the case where $\tau=1$. 
			This is an example of the Artin-Mazur zeta function \cite{artinmazur65} 
			corresponding to the dynamical system naturally defined on a finite simple graph. 
			We call this original one the {\it Bowen-Lanford zeta function}. 
			In the case where $\tau=1$, the map $\theta^{\rm G}$ is denoted by $\theta^{\rm BL}$(=$\theta^{\rm BL}_1$). 
			We can see that the map $\theta^{{\rm BL}}_{\tau}$ satisfies the adjacency condition. 
			Thus $(\Pi_{\Delta}, \theta^{\rm BL}_\tau)$ satisfy the path condition by Lemma \ref{lem : adjacency conditions}. 
			Hence we have the identities 
			$
				Z_{\Delta}(t; \theta^{{\rm BL}}_{\tau})
				=
				E_{\Delta}(t; \theta^{{\rm BL}}_{\tau})
				=
				H_{\Delta}(t; \theta^{{\rm BL}}_{\tau}).
			$
			Note also that the pair $(\Pi_{\Delta}^{\flat}, \theta^{{\rm BL}}_{\tau})$ 
			does not satisfy the path condition. 
			Hence, although the identity 
			$
				Z^{\flat}_{\Delta}(t; \theta^{{\rm BL}}_{\tau})
				=
				E^{\flat}_{\Delta}(t; \theta^{{\rm BL}}_{\tau})
			$ holds, 
			the formal power series 
			$
				H^{\flat}_{\Delta}(t;\theta^{{\rm BL}}_{\tau})
			$ 
			is not equal to these. 
			See \cite{morita20} for precise information. 
			(Remark that, in \cite{morita20}, 
			the weighed Bowen-Lanford zeta function is called the Mizuno-Sato zeta function.)

			It is also known that the Bowen-Lanford zeta $Z_{\Delta}(t; \theta^{{\rm BL}})$ has the Ihara expression. 
			Let $\Gamma$ be a finite simple graph, 
			and $A_{\Gamma}$ the adjacency matrix of $\Gamma$. 
			We denote the following formal power series 
			$$
				\frac{1}{\det(I-tA_{\Gamma})}.
			$$ 
			by $I_{\Gamma}(t; \theta^{{\rm BL}})$. 
			In \cite{bowenlanford70}, 
			R. Bowen and O. Lanford  shows 
			the identity $Z_{\Gamma}(t; \theta^{{\rm BL}})=I_{\Gamma}(t; \theta^{{\rm BL}})$, 
			which we call the Ihara expression of the Bowen-Lanford zeta $Z_{\Gamma}(t; \theta^{{\rm BL}})$. 
			 See also \cite{mizunosato01}.

			\subsubsection{The Mizuno-Sato zeta function}
			\label{subsubsection : The Mizuno-Sato zeta function}
			In the case where $\tau=\upsilon$, 
			the generalized weighted zeta $Z_{\Delta}(t; \theta^{\rm G})$ is called 
			the {\it Mizuno-Sato zeta function} of $\Delta$, 
			which was originally introduced in \cite{mizunosato03a}. 
			In this case,
			we denote the map $\theta^{\rm G}$ by $\theta^{\rm MS}$. 
			The map $\theta^{\rm MS}$ satisfies the reduced adjacency condition 
			by virtue of Lemma \ref{lem : The adjacency condition for G}. 
			Thus both $(\Pi_{\Delta}, \theta^{\rm MS})$ and $(\Pi^{\flat}_{\Delta}, \theta^{\rm MS})$ 
			satisfy the path condition by Lemma \ref{lem : adjacency conditions}. 
			Hence we have the identities 
			$
				Z_{\Delta}(t; \theta^{\rm MS})
				=
				E_{\Delta}(t; \theta^{\rm MS})
				=
				H_{\Delta}(t; \theta^{\rm MS}), 
			$
			and 
			$
				Z^{\flat}_{\Delta}(t; \theta^{\rm MS})
				=
				E^{\flat}_{\Delta}(t; \theta^{\rm MS})
				=
				H^{\flat}_{\Delta}(t; \theta^{\rm MS}).
			$
			Since $\theta^{\rm MS}$ excludes the non-reduced closed paths, 
			we have $Z_{\Delta}(t; \theta^{\rm MS})=Z^{\flat}_{\Delta}(t; \theta^{\rm MS})$, 
			and those six formal power series are all equal to each other. 
			Note that, for both cases $\Pi_{\Delta}$ and $\Pi^{\flat}_{\Delta}$, 
			we have $N_m(\theta^{\rm MS})=\sum_{x\in X_m^{\flat}}\cir_{\theta^{\rm MS}}(x)$. 
			See 4.3 in \cite{morita20}.

			The Ihara expression of the Mizuno-Sato zeta is given as follows \cite{mizunosato04}. 
			Let $\Gamma=(V, E)$ be a finite simple graph. 
			Note that, since $\Gamma$ is simple, there exists at most one edge $e$ 
			lying between two distinct vertices $u, v\in V$. 
			We denote this possible single edge $e$ by $\{u, v\}$. 
			Let $\Delta(\Gamma)=(V, {\cal A})$ be the symmetric digraph of $\Delta$. 
			The arc set ${\cal A}$ consists of arcs $a_e, \overline{a_e}$ for $e\in E$. 
			Given a total order $<$ on $V$, 
			if $e=\{u, v\}$ with $u<v$, 
			then we say $a_e\in{\cal A}_{uv}$ and denote it by $a_{uv}$, 
			$\overline{a_e}$ by $a_{vu}$. 
			Consider the matrix $W_{\Gamma}=(\tau(a_{uv}))_{u, v\in V}$, 
			called the {\it weighted matrix} of $\Gamma$ \cite{mizunosato04}. 
			Here we understand that $\tau(a_{uv})=\tau(a_{vu})=0$ if $\{u, v\}\notin E$. 
			Let $I_{\Gamma}(t; \theta^{\rm MS})$ denote the formal power series 
			given by the reciprocal of the following polynomial 
			$$
				(1-t^2)^{|E|-|V|}
				\det(I-tW_{\Gamma}+t^2(D_{\Gamma}-I)).
			$$
			In \cite{mizunosato04}, H. Mizuno and I. Sato shows 
			the identity $Z_{\Gamma}(t; \theta^{\rm MS})=I_{\Gamma}(t; \theta^{\rm MS})$, 
			the Ihara expression of $Z_{\Gamma}(t; \theta^{\rm MS})$.

			\subsubsection{The Sato zeta function}
			\label{subsubsection : The Sato zeta function}
			In the case where $\upsilon=1$, 
			the generalized weighted zeta $Z_{\Delta}(t; \theta^{\rm G})$ is called 
			the {\it Sato zeta function}\cite{sato07} (see also \cite{morita20}, 4.4). 
			In this case,
			we denote the map $\theta^{\rm G}$ by $\theta^{\rm S}$. 
			The map $\theta^{\rm S}$ satisfies the adjacency condition. 
			Thus $(\Pi_{\Delta}, \theta^{\rm S})$ 
			satisfy the path condition. 
			Hence we have the identities 
			$
				Z_{\Delta}(t; \theta^{\rm S})
				=
				E_{\Delta}(t; \theta^{\rm S})
				=
				H_{\Delta}(t; \theta^{\rm S}).
			$
			Note that the map $\theta^{\rm S}$ does not satisfy 
			the reduced adjacency condition in general (c.f., Lemma \ref{lem : The adjacency condition for G}). 
			In this case, 
			thought the identity 
			$
				Z^{\flat}_{\Delta}(t; \theta^{\rm S})
				=
				E^{\flat}_{\Delta}(t; \theta^{\rm S})
			$ 
			holds, 
			these are not equal to $H^{\flat}_{\Delta}(t; \theta^{\rm S})$. 
			Let $\Gamma=(V, E)$ be a finite simple graph, 
			and $\Delta=(V, {\cal A})$ the symmetric digraph of $\Gamma$. 
			Let $W_{\Gamma}=(w_{uv})_{u, v\in V}$ denote the weighted matrix of $\Gamma$ 
			as in \ref{subsubsection : The Mizuno-Sato zeta function}, 
			and 
			$
				D_{\Gamma}(w)
				=
				(d_{uu'})_{u, u'v\in V}
			$ 
			the diagonal matrix given by 
			$
				d_{uu'}
				=
				\delta_{uu'}
				\sum_{a_{uv}\in{\cal A}_{u*}}
				w_{uv}.
			$ 
			If $I_{\Delta}(t; \theta^{\rm S})$ denote the reciprocal of the polynomial
			$$
				(1-t^2)^{|E|-|V|}
				\det(I-tW_{\Gamma}+t^2(D_{\Gamma}(w)-I)),
			$$
			the Ihara expression of the Sato zeta $Z_{\Gamma}(t; \theta^{\rm S})$ 
			is given by the identity $Z_{\Gamma}(t; \theta^{\rm S})=I_{\Gamma}(t; \theta^{\rm S})$ \cite{sato07}.

			The Sato zeta function $Z_{\Delta}(t; \theta^{\rm S})$ 
			have a relation with quantum walks on finite graphs \cite{konnosato12}. 
			Let $\Gamma$ be a finite simple graph, 
			and $\Delta=\Delta(\Gamma)$ the symmetric digraph of $\Gamma$. 
			If we let $\tau(a')=1/\deg {\frak t}(a')$, 
			the Sato zeta $Z_{\Gamma}(t; \theta^{\rm S})$ gives 
			the reciprocal characteristic polynomial of 
			the time evolution matrix $U$ of 
			the \lq\lq Grover walk\rq\rq on $\Gamma$ \cite{grover96}, 
			a quantum walk model defined on a finite simple graph $\Gamma$ 
			which is extensively studied in the area. 
			Konno and Sato \cite{konnosato12} shows that, 
			motivated by Emms et.al.\cite{EHSW06a}, 
			the Ihara expression for the Sato zeta $Z_{\Delta}(t; \theta^{\rm S})$ (see \cite{sato07}) 
			gives a fine description of the spectrum of $U$. 
			By the Konno-Sato theorem\cite{konnosato12}, 
			we may understand that 
			the Hashimoto expression of a graph zeta for $\Gamma$ depicts the time evolution 
			of a quantum walk model on $\Gamma$, 
			and the Ihara expression describes the spectral mapping theorem for the model. 
			For the spectral mapping theorem, 
			consult \cite{matueogurisusegawa17, segawasuzuki19} (see also \cite{ishikawa21+}). 
			In particular, the result suggests that 
			the edge matrix $M_{\Delta(\Gamma)}(\theta^{\rm G})$ may provide 
			the time evolution matrix of a quantum walk model on $\Gamma$, 
			and may produce a family of quantum walks on $\Gamma$ 
			which contains the Grover walk as an example. 
			A. Ishikawa \cite{ishikawa21+} considers this problem and 
			obtains such a family of quantum walks on finite graphs.

			\subsubsection{The Bartholdi zeta function}
			\label{subsubsection : The Bartholdi zeta function}

			Let $q$ be an indeterminate, 
			and consider the polynomial algebra $R[q]$, 
			a commutative ${\mathbb Q}$-algebra with unity. 
			Let $\theta^{\rm B}$ be a map ${\cal A}\times{\cal A}\rightarrow R[q]$ 
			given by 
			$$
			\theta^{\rm B}(a, a')
			=
			\delta_{\head(a)\tail(a')}-(1-q)\delta_{a'\in S(a)}.
			$$ 
			This is the case where $\tau(a)=1$ and $\upsilon(a)=1-q$ ($a\in{\cal A}$) in $\theta^{\rm G}$. 			
			The generalized weighted zeta $Z_{\Delta}(t; \theta^{\rm B})$ is called 
			the {\it Bartholdi zeta function} \cite{bartholdi99} of $\Delta$ (see also \cite{mizunosato03a}). 
			The map $\theta^{\rm B}$ satisfies the adjacency condition. 
			Thus $(\Pi_{\Delta}, \theta^{\rm B})$ satisfies the path condition. 
			Hence we have the identities 
			$
				Z_{\Delta}(t; \theta^{\rm B})
				=
				E_{\Delta}(t; \theta^{\rm B})
				=
				H_{\Delta}(t; \theta^{\rm B}).
			$ 
			On $\Pi_{\Delta}^{\flat}$, 
			though the identity 
			$
				Z^{\flat}_{\Delta}(t; \theta^{\rm B})
				=
				E^{\flat}_{\Delta}(t; \theta^{\rm B})
			$ 
			holds, 
			these does not equals $H^{\flat}_{\Delta}(t; \theta^{\rm B})$ 
			since the map $\theta^{\rm B}$ does not satisfy 
			the reduced adjacency condition.

			The Euler expression $E_{\Delta}(t; \theta^{\rm B})$ is described by 
			a combinatorial statistics, called 
			the \lq\lq cyclic bump count\rq\rq. 
			Let $x=(x_i)_{i\in{\mathbb Z}}\in X_m$ be an $m$-periodic point in $\Pi_{\Delta}$, 
			and let ${\rm cbc}(x)$ denote the cardinality of the set 
			$$
				\{i=k, k+1, \dots, k+m-1\mid x_{i+1}\in S(x_i)\},
			$$ 
			where $(x_k, x_{k+1}, \dots, x_{k+m-1})$ is an $m$-section of $x$. 
			Note that $\cbc(x)$ does not depend on the choice of an $m$-section for $x\in X_m$. 
			If $a'\in S(a)$, 
			then we have $\theta^{\rm B}(a, a')=q$. 
			This lead to the identity 
			$
				\cir_{\theta^{\rm B}}(x)
				=
				q^{\cbc(x)}
			$ 
			for each $x\in X_m$, 
			and it follows from Proposition \ref{prop : The Euler expression for graph zetas} 
			that $E_{\Delta}(t; \theta^{\rm B})$ equals 
			$$
				\prod_{[x]\in[X]}
				\frac{1}{1-q^{\cbc(\pi(x))}t^{\varpi(x)}}.
			$$ 
			By Proposition \ref{prop : The Hashimoto expression for graph zetas}, 
			the Hashimoto expression is given by 
			$H_{\Delta}(t; \theta^{\rm B})=1/\det(I-tM_{\Delta}(\theta^{\rm B}))$, 
			where $M_{\Delta}(\theta^{\rm B})=(\theta^{\rm B}(a, a'))_{a, a'\in{\cal A}}$. 
			The Bartholdi zeta was originally considered for the symmetric digraph 
			of a finite graph, 
			and it is defined by the Euler expression. 
			See \cite{bartholdi99} for the original definition. 
			Let $\Gamma=(V, E)$ be a finite graph, 
			and let $L$ denote the subset of $E$ consisting of loops in $\Gamma$. 
			The Ihara expression of the Bartholdi zeta is the identity 
			$
				Z_{\Gamma}(t; \theta^{\rm B})
				=
				I_{\Gamma}(t; \theta^{\rm B}), 
			$ 	
			where $I_{\Delta}(t; \theta^{\rm B})\in R[[t]]$ is the reciprocal of the following polynomial 
			$$
				(1+(1-q)t)^{|L|}
				(1-(1-q)^2t^2)^{|E|-|L|-|V|}
				\det(I-tA_{\Gamma}+(1-q)t^2(D_{\Gamma}-(1-q)I)).
			$$
			See \cite{bartholdi99} for precise information.

			Since 
			\begin{equation}\label{equation : Bartholdi map at q=0}
				\left.
					\theta^{\rm B}
				\right|_{q=0}
				=
				\theta^{\rm I}, 
			\end{equation}
			the Bartholdi zeta $Z_{\Delta}(t; \theta^{\rm B})$ equals 
			the Ihara zeta $Z_{\Delta}(t; \theta^{\rm I})$ at $q=0$. 
			In this sense, we also call $Z_{\Delta}(t; \theta^{\rm B})$ 
			the {\it Ihara zeta function} of {\it Bartholdi type}, 
			and denote it by $Z_{\Delta}(q, t; \theta^{\rm I})$, 
			the map $\theta^{\rm B}$ by $\theta^{\rm I}_q$. 
			Hence we have $Z_{\Delta}(q, t; \theta^{\rm I})=Z_{\Delta}(t; \theta^{\rm I}_q)$. 
			In the same sense, the Euler expression $E_{\Delta}(t; \theta^{\rm I}_q)$ and 
			the Hashimoto expression $H_{\Delta}(t; \theta^{\rm I}_q)$ 
			are denoted by $E_{\Delta}(q, t; \theta^{\rm I})$ and $H_{\Delta}(q, t; \theta^{\rm I})$, 
			respectively. 
			If we agree that $0^0=1$, 
			then we have
			\begin{equation}\label{equation : cyclic bump count at q=0}
				q^{\cbc(x)}
				=
				\left\{
					\begin{array}{l}
					1,\quad\mbox{if $x\in X_m^{\flat}$}, \\
					0, \quad\mbox{if $x\in X_m\setminus X_m^{\flat}$}.
					\end{array}
				\right.
			\end{equation}
			This leads to the identity 
			$
				\sum_{x\in X_m}\cir_{\theta^{\rm I}_q|_{q=0}}(x)
				=\sum_{x\in X_m^{\flat}}\cir_{\theta^{\rm I}}(x), 
			$
			i.e., 
			$
				N_m(\theta^{\rm I}_q|_{q=0})
				=
				N_m(\theta^{\rm I}), 
			$
			which directly shows that $Z_{\Delta}(0, t; \theta^{\rm I})=Z_{\Delta}(t; \theta^{\rm I})$. 
			The identity $E_{\Delta}(0, t; \theta^{\rm I})=E_{\Delta}(t; \theta^{\rm I})$ 
			also follows directly from (\ref{equation : cyclic bump count at q=0}) . 
			The identity $H_{\Delta}(0, t; \theta^{\rm I})=H_{\Delta}(t; \theta^{\rm I})$ 
			follows trivially from (\ref{equation : Bartholdi map at q=0}). 
			If we let $q=1$, 
			then $Z_{\Delta}(q, t; \theta^{\rm I})$ equals the Bowen-Lanford zeta $Z_{\Delta}(t; \theta^{\rm BL})$. 
			Thus the Bartholdi zeta interpolates the Ihara zeta and the Bowen-Lanford zeta.

			\vspace*{5mm}
			There exists other classical graph zetas which were not discussed here. 
			For instance, it should be mentioned here about the {\it edge zeta function} \cite{starkterras96} 
			and the {\it path zeta function} (see Terras \cite{terras11}). 
			These are also graph zetas in the sense of Definition \ref{def : Graph zeta functions}. 
			Therefore they have the three expression on the fundamental framework in our development. 
			See 4.5 of \cite{morita20} for further information.

		\subsection{Graph zeta functions of Bartholdi type}
		\label{subsection : Graph zeta functions of Bartholdi type}
		
		In \ref{subsubsection : The Bartholdi zeta function}, 
		we introduce the Ihara zeta of Bartholdi type. 
		In the same manner, 
		we can define the Bartholdi type for various graph zetas. 
		Let $\Delta=(V, {\cal A})$, $S$, $R$, $\tau$, $\upsilon$ and $\theta^{\rm G}$ be as in 
		\ref{subsection :　Definition and the three expressions}. 
		Let $q$ be an indeterminate and 
		$$
			\theta^{\rm G}_q : {\cal A}\times {\cal A}\rightarrow R[q]
		$$
		a map given by 
		$
			\theta^{\rm G}_q(a, a')
			=
			\tau(a')\delta_{\head(a)\tail(a')}
			-
			(1-q)\upsilon(a')\delta_{a'\in S(a)}.
		$
		Let $\Xi$ be a $\varphi$-stable subset of ${\cal A}^{\mathbb Z}$
		
		\begin{df}[Bartholdi type]
		\label{def : Bartholdi type}
		{\em
		The formal power series $Z_{\Xi}(t; \theta^{\rm G}_q)$ is called 
		the generalized weighted zeta function of {\it Bartholdi type}, 
		which is also denoted by $Z_{\Xi}(q, t; \theta^{\rm G})$. 
		The Euler expression $E_{\Xi}(t; \theta^{\rm G}_q)$ 
		is denoted by $E_{\Xi}(q, t; \theta^{\rm G})$, 
		and the Hashimoto expression $H_{\Xi}(t; \theta^{\rm G}_q)$ 
		by $H_{\Xi}(q, t; \theta^{\rm G})$ if exists. 
		}
		\end{df}
		
		\begin{lem}
		\label{lem : The adjacency condition for Gq}
		The map $\theta^{\rm G}_q$ satisfies the adjacency condition. 
		\end{lem}
		{\it Proof.}
		Suppose that $\theta^{\rm G}_q(a, a')\neq 0$. 
		If $a'\in S(a)$, then we obviously have $\head(a)=\tail(a')$. 
		If $a'\notin S(a)$, then we have 
		$
			\theta^{\rm G}_q(a, a')
			=
			\tau(a')\delta_{\head(a)\tail(a')}, 
		$
		and the assumption forces $\head(a)=\tail(a')$. 
		\qed
		\begin{prop}
		\label{prop : The three expressions for the generalized weighted zeta function of Bartholdi type}
		We have the identities 
		$
			Z_{\Delta}(q, t; \theta^{\rm G})
			=
			E_{\Delta}(q, t; \theta^{\rm G})
			=
			H_{\Delta}(q, t; \theta^{\rm G}).
		$
		\end{prop}
		{\it Proof.}
		By Proposition \ref{prop : The Euler expression for graph zetas}, 
		it follows that 
		$
			Z_{\Delta}(q, t; \theta^{\rm G})
			=
			E_{\Delta}(q, t; \theta^{\rm G}).
		$ 
		By Proposition \ref{prop : The Hashimoto expression for graph zetas}, 
		it follows that 
		$
			E_{\Delta}(q, t; \theta^{\rm G})
			=
			H_{\Delta}(q, t; \theta^{\rm G}), 
		$ 
		since $(\Pi_\Delta, \theta^{\rm G}_q)$ satisfies the path condition 
		by Lemma \ref{lem : The adjacency condition for Gq}. 
		\qed
		
		\vspace*{5mm}
		Varying $\tau$ and $\upsilon$, 
		we obtain various graph zetas of Bartholdi type. 
		In the case where $\tau=\upsilon$ for example, 
		we have 
		$
			\left.\theta^{\rm G}_q\right|_{\tau=\upsilon}(a, a')
			=
			\tau(a')(\delta_ {\head(a)\tail(a')}-(1-q)\delta_{a'\in S(a)}). 
		$ 
		If $q=0$, then we have 
		$$
			\left.\theta^{\rm G}_q\right|_{\tau=\upsilon\atop q=0}
			=
			\theta^{\rm MS}.
		$$
		Thus we denote the map $\left.\theta^{\rm G}_q\right|_{\tau=\upsilon}$ 
		by $\theta^{\rm MS}_q$ for a generic $q$, 
		and the graph zeta $Z_{\Delta}(t; \theta^{\rm MS}_q)$ 
		by $Z_{\Delta}(q, t; \theta^{\rm MS})$, 
		which we call the {\it Mizuno-Sato zeta function} of {\it Bartholdi type}
		(c.f., \ref{subsubsection : The Bartholdi zeta function}.) 
		This zeta appears in \cite{choekwakparksato07} 
		for the case where $\Delta$ is simple. 
		We also denote the Euler expression $E_{\Delta}(t; \theta^{\rm MS}_q)$ 
		and the Hashimoto expression $H_{\Delta}(t; \theta^{\rm MS}_q)$ by 
		$E_{\Delta}(q, t; \theta^{\rm MS})$ and $H_{\Delta}(q, t; \theta^{\rm MS})$, 
		respectively. 
		Let $x=(a_i)_{i\in{\mathbb Z}}\in X_m$ be an $m$-periodic point. 
		If $a_{i+1}\in S(a_i)$, 
		then we have $\theta^{\rm G}_q(a_i, a_{i+1})=\tau(a_{i+1})q$; 
		otherwise, 
		$\theta^{\rm G}_q(a_i, a_{i+1})=\tau(a_{i+1})$. 
		Hence we have 
		\begin{equation}\label{equation : The circular product for Gq}
			\cir_{\theta^{\rm MS}_q}(x)
			=
			\cir_{\theta^{\rm BL}_{\tau}}(x)q^{\cbc(x)},
		\end{equation} 
		and 
		$
			N_m(\theta^{\rm MS}_q)
			=
			\sum_{x\in X_m}
			\cir_{\theta^{\rm BL}_{\tau}}(x)q^{\cbc (x)},
		$
		which equals $N_m(\theta^{\rm BL}_{\tau})$ if $q=1$. 
		The agreement $0^0=1$ (c.f., (\ref{equation : cyclic bump count at q=0})) leads to the identity 
		$
			N_m(\left.\theta^{\rm MS}_q\right|_{q=0})
			=
			N_m(\theta^{\rm MS}).
		$
		Thus the exponential expression $Z_{\Delta}(q, t; \theta^{\rm MS})$ at $q=0$ and $q=1$ 
		are actually identical with 
		$Z_{\Delta}(t; \theta^{\rm MS})$ and $Z_{\Delta}(t; \theta^{\rm BL}_\tau)$ respectively. 
		It follows from Proposition \ref{prop : The Euler expression for graph zetas} 
		that the Euler expression $E_{\Delta}(q, t; \theta^{\rm MS})$ is given by 
		$$
			\prod_{[x]\in[X]}
			\frac{1}{1-q^{\cbc(x)}\cir_{\theta^{\rm BL}_{\tau}(x)}t^{\varpi(x)}}. 
		$$
		By Proposition \ref{prop : The Hashimoto expression for graph zetas}, 
		the Hashimoto expression $H_{\Delta}(q, t; \theta^{\rm MS})$ is given by 
		$1/\det(I-tM)$, 
		where $M=(\theta^{\rm MS}_q(a, a'))_{a, a'\in{\cal A}}$.

		On the other hand, 
		the map $\theta^{\rm G}_q|_{\upsilon =1}$ defines 
		another graph zeta of Bartholdi type. 
		In the case of $\upsilon=1$, we denote $\theta^{\rm G}_q$ by $\theta^{\rm S}_q$, 
		and the graph zeta $Z_{\Delta}(t; \theta^{\rm S}_q)$ is called 
		the {\it Sato zeta function} of {\it Bartholdi type}. 
		We also denote $Z_{\Delta}(t; \theta^{\rm S}_q)$ by $Z_{\Delta}(q, t; \theta^{\rm S})$. 
		By Proposition \ref{prop : The Euler expression for graph zetas} 
		and Proposition \ref{prop : The Hashimoto expression for graph zetas}, 
		the exponential expression $Z_{\Delta}(q, t; \theta^{\rm MS})$ can be reformulated into 
		the Euler expression $E_{\Delta}(q, t; \theta^{\rm S})$ 
		and the Hashimoto expression $H_{\Delta}(q, t; \theta^{\rm S})$ 
		in the same manner.

		To summarize, 
		as is suggested by these examples, 
		if we vary $\tau$ and $\upsilon$ in $\theta^{\rm G}_q$, 
		then we obtain various graph zetas of Bartholdi type. 
		In particular, if we consider the case $\Xi=\Pi_{\Delta}$, 
		it follows from Proposition 
		\ref{prop : The three expressions for the generalized weighted zeta function of Bartholdi type} 
		these zetas of Bartholdi type are always combinatorial, 
		that is, they always have the Hashimoto expression.

	\section{The Ihara expression for the generalized weighted zeta function of Bartholdi type}
		
		In this section, 
		the Ihara expression of the generalized weighted zeta of Bartholdi type is given. 
		Since $\theta^{\rm G}_q=\theta^{\rm G}$ at $q=0$, 
		the result for the case of the generalized weighted zeta is included in the main theorem. 
		However, it is enough to show the main theorem for $\theta=\theta^{\rm G}$ 
		since the polynomial ring $R[q]$ is also a commutative ${\mathbb Q}$-algebra with unity. 
		If we choose $R=R[q]$ and replace $\upsilon$ by $(1-q)\upsilon$ 
		in the identity (\ref{equation : The generalized weighted map}), 
		then the map $\theta^{\rm G}_q$ is obtained. 
		In the following of this section, 
		let $\Delta=(V, {\cal A})$ be a finite digraph, 
		$R$ a commutative ${\mathbb Q}$-algebra with unity, 
		and $\theta=\theta^{\rm G}$. 
		Note that the digraph $\Delta$ allows multi-arcs and multi-loops. 
		Since $\theta^{\rm G}$ satisfies the adjacency condition (Lemma \ref{lem : The adjacency condition for G}), 
		the generalized weighted zeta $Z_{\Delta}(t; \theta^{\rm G})$ has 
		the Hashimoto expression $H_{\Delta}(t; \theta^{\rm G})$ 
		(Corollary \ref{cor : The three expression for the generalized weighted zeta}). 
		Recall that 
		$
			H_{\Delta}(t; \theta^{\rm G})
			=
			1/\det(I-tM_{\Delta}(\theta^{\rm G}))
		$, 
		where 
		$
			M_{\Delta}(\theta^{\rm G})
			=
			(\theta^{\rm G}(a, a'))_{a, a'\in{\cal A}}.
		$

		\subsection{Auxiliary results}
		
		\subsubsection{The Schur complement}
		\label{subsubsection : The Schur complement}
		
		We recall some auxiliary facts on block matrices. 
		Let 
		\begin{equation}\label{equation : The block matrix}
			M=
			\left[
				\begin{array}{cc}
				A & B\\
				C & D
				\end{array}
			\right]
		\end{equation}
		be a block matrix, 
		where $A$ and $D$ are square matrices. 
		Suppose that $A$ and $D$ are invertible. 
		The results in \ref{subsubsection : The Schur complement} are well-known. 
		However, brief proofs or comments are provided for readers convenience.

		\begin{df}[The Schur complement]
		\label{def : The Schur complement}
		{\em 
		Let 
		$$
			\begin{array}{l}
			M/A:=D-CA^{-1}B,\\
			M/D:=A-BD^{-1}C.
			\end{array}
		$$
		The matrix $M/A$ (resp. $M/D$) is called 
		the {\it Schur complement} of $A$ (resp. $D$) in $M$.
		}
		\end{df}

		\begin{lem}\label{lem : The Woodbury identity}
		The Schur complement $M/A$ is invertible if and only if $M/D$ is invertible. 
		In this case, we have 
		\begin{equation}\label{equation : The Woodbury identity for A}
			(M/A)^{-1}=D^{-1}+D^{-1}C(M/D)^{-1}BD^{-1},
		\end{equation}
		\begin{equation}\label{equation : The Woodbury identity for D}
			(M/D)^{-1}=A^{-1}+A^{-1}B(M/A)^{-1}CA^{-1}.
		\end{equation}
		\end{lem}

		It is possible to see that $(M/A)\{D^{-1}+D^{-1}C(M/D)^{-1}BD^{-1}\}=I$ and vise versa. 
		The identity (\ref{equation : The Woodbury identity for D}) is also similarly verified. 
		These identities (\ref{equation : The Woodbury identity for A}),  (\ref{equation : The Woodbury identity for D}) 
		are called the {\it Woodbury identity}. 
		Carrying out the multiplication of block matrices, 
		one can easily show the following lemma.

		\vspace*{5mm}
		\begin{lem}\label{lem : The inverse for a block matrix}
		If $A$ and $D$ are invertible in the block matrix $M$, 
		then we have the following two decompositions:
		\begin{eqnarray}
		M 
		&=&
			\left[
				\begin{array}{cc}
				I & O\\CA^{-1} & I
				\end{array}
			\right]
			\left[
				\begin{array}{cc}
				A & O\\O & M/A
				\end{array}
			\right]
			\left[
				\begin{array}{cc}
				I & A^{-1}B\\O & I
				\end{array}
			\right]
			\label{equation : The determinant of the block matrix 1}
			\\
		&=&
			\left[
				\begin{array}{cc}
				I & BD^{-1}\\O & I
				\end{array}
			\right]
			\left[
				\begin{array}{cc}
				M/D & O\\O & D
				\end{array}
			\right]
			\left[
				\begin{array}{cc}
				I & O\\D^{-1}C & I
				\end{array}
			\right].
			\label{equation : The determinant of the block matrix 2}
		\end{eqnarray}
		In particular, 
		we have 
		$
			\det M=\det A\det M/A=\det M/D \det D.
		$
		\end{lem}

		\vspace*{5mm}
		If $A$ and $M/A$ are invertible, 
		the inverse of $M$ is described as follows. 
		In the case where $D$ and $M/D$ are invertible, 
		one can obtain a similar formula. 
		Recall that 
		\begin{equation}\label{equation : The inverse of triangular block matrices}
			\left[
			\begin{array}{cc}
				I & B\\
				O & I
			\end{array}
			\right]^{-1}
			=
			\left[
			\begin{array}{cc}
				I & -B\\
				O & I
			\end{array}
			\right], \quad
			\left[
			\begin{array}{cc}
				I & O\\
				C & I
			\end{array}
			\right]^{-1}
			=
			\left[
			\begin{array}{cc}
				I & O\\
				-C & I
			\end{array}
			\right]. 
		\end{equation}

		\begin{cor}\label{cor : The inverse for the block matrix}
		If $A$ and $M/A$ are invertible, 
		then $M$ is invertible and 
		the inverse of the block matrix $M$ is given by 
		$$
			M^{-1}
			=
			\left[
				\begin{array}{cc}
				(M/D)^{-1} & -A^{-1}B(M/A)^{-1}\\
				-(M/A)^{-1}CA^{-1} & (M/A)^{-1}
				\end{array}
			\right]. 
		$$
		\end{cor}
		{\it Proof.}
		By Lemma \ref{lem : The inverse for a block matrix} and 
		the identities (\ref{equation : The inverse of triangular block matrices}), 
		we have 
		$$
			M^{-1}=
			\left[
				\begin{array}{cc}
				I & -A^{-1}B\\
				O & I
				\end{array}
			\right]
			\left[
				\begin{array}{cc}
				A^{-1} & O\\
				O & (M/A)^{-1}
				\end{array}
			\right]
			\left[
				\begin{array}{cc}
				I & O\\
				-CA^{-1} & I
				\end{array}
			\right].
		$$
		The assertion follows from Lemma \ref{lem : The Woodbury identity}.\qed

		\vspace*{5mm}
		
		In the proof of the main theorem, 
		we concentrate on a block matrix of the form 
		\begin{equation}\label{equation : The block matrix of the particular form}
			M=
			\left[
				\begin{array}{cc}
				I & B\\
				C & I
				\end{array}
			\right], 
		\end{equation} 
		where $A$ and $D$ are identity matrices in the block matrix (\ref{equation : The block matrix}). 
		In this case, the Schur complements are given by 
		$$
			M/A=I-CB, 
			\quad 
			M/D=I-BC.
		$$

		\begin{cor}\label{cor : The inverse of the block matrix with Schur complements}
			If $M$ is a matrix of the form {\em (\ref{equation : The block matrix of the particular form})}, 
			then $I-BC$ is invertible if and only if $I-CB$ is invertible. 
			In this case, 
			the matrix $M$ itself is also invertible, 
			and the inverse $M^{-1}$ is given by 
			\begin{equation}
				\left[
					\begin{array}{cc}
					(I-BC)^{-1} & -B(I-CB)^{-1}\\
					-(I-CB)^{-1}C & (I-CB)^{-1}
					\end{array}
				\right].
			\end{equation}
			We also have 
			$
				\det M
				=
				\det(I-BC)
				=
				\det(I-CB).
			$ 
			
		\end{cor}
		{\it Proof.}
		The first assertion follows from Lemma \ref{lem : The Woodbury identity}. 
		The second assertion follows from Corollary \ref{cor : The inverse for the block matrix}. 
		The third assertion follows from Lemma \ref{lem : The inverse for a block matrix}. 
		\qed

		\subsubsection{Column constant matrices}
		\label{subsubsection : Column constant matrices}
		If a matrix $M=(m_{ij})$ satisfies the condition that 
		$
			m_{rj}=m_{sj}
		$ 
		for all $r, s, j$, 
		then $M$ is called a {\it column constant matrix}. 
		For a column constant matrix $M=(m_{ij})$, 
		the sum $\sum_{j\geq 1}m_{ij}$ does not depend on $i$, 
		which we call the {\it row sum} of $M$. 
		We denote the row sum of $M$ by $\rho_M$. 
		It is clear that the row sum $\rho_M$ of a column constant matrix $M$ 
		equals the trace ${\rm tr}\,M$ of $M$.

		\begin{lem}\label{lem : The square of a column constant matrix}
		If a square matrix $M$ is column constant, 
		then we have $M^2=\rho_MM$. 
		\end{lem}
		{\it Proof.}
		It follows from a direct calculation.\qed

		\begin{lem}\label{lem : The product of column constant matrices}
		Let $M=(m_{ij})$, $N=(n_{ij})$ be matrices with sizes $k\times l$ and $l\times k$ respectively. 
		If $M$ and $N$ are column constant, then 
		the products $MN$, $NM$ are also column constant. 
		If $MN$, $NM$ are square matrices, 
		then ${\rm tr}\,(MN)$ and ${\rm tr}\,(NM)$ equal 
		the product $\rho_{MN}=\rho_M\rho_N$ 
		of the row sums of $M$ and $N$. 
		\end{lem}
		{\it Proof.}
			It suffices to show the assertion for $MN$. 
			Since $N$ is column constant, 
			the entries $n_{ij}$ belonging to the $j$-th column of $N$ does not depend on the indices $i$, 
			and we denote it by $n_j$, i.e., $n_j=n_{ij}$ for all $i$. 
			The $(i, j)$-entry of $MN$ is given by 
			\begin{equation}\label{equation : The (i, j)-entries of the product of two column constant matrices}
				\sum_{k\geq 1}m_{ik}n_{kj}
				=
				\rho_Mn_j, 
			\end{equation}
			which does not depend on $i$. 
			Thus the product $MN$ is column constant matrix. 
			The identity (\ref{equation : The (i, j)-entries of the product of two column constant matrices}) 
			also shows that 
			$
				{\rm tr}\,(MN)
				=
				\sum_{i=1}^k
				\rho_Mn_i
			$, 
			which equals $\rho_M\rho_N$. 
			\qed

		\vspace*{5mm}
		If a matrix $M=(m_{ij})\in M_n(R)$ is column constant, 
		then one can easily see that 
		$\det(I+M)=1+(m_{11}+m_{22}+\cdots +m_{nn})$.

		\begin{lem}\label{lem : The determinant of a column constant matrix}
			We have $\det(I+M)=1+\rho_M$ 
			for any column constant square matrix $M$. 
		\end{lem}

		Let $t$ be an indeterminate, 
		and $R[[t]]$ the ring of formal power series with coefficients in $R$. 
		For a column constant square matrix $M$, 
		by Lemma \ref{lem : The determinant of a column constant matrix}, 
		we have 
		$
			\det(I+tM)=1+t\rho_M, 
		$ 
		which is an invertible element in $R[[t]]$. 
		
		\begin{lem}\label{lem : The inverse of a column constant matrix}
		Let $M\in M_n(R)$. 
		If $M$ is column constant, 
		then $I+tM$ is invertible as an element of $M_n(R[[t]])$, 
		and the inverse matrix $(I+tM)^{-1}$ is given by 
		$$
			\det(I+tM)^{-1}
			\left\{I-t(M-\rho_MI)\right\}.
		$$
		\end{lem}
		{\it Proof.}
			If we multiply the matrix $I+tM$ by $\det(I+tM)^{-1}(I-t(M-\rho_MI)$ from the left, 
			then it equals 
			$$
			(1+t\rho_M)^{-1}
			(I+t\rho_MI-t^2M^2+t^2\rho_MM).
			$$
			By Lemma \ref{lem : The square of a column constant matrix}, 
			this equals 
			$
			(1+t\rho_M)^{-1}
			(I+t\rho_MI)=I. 
			$ 
			The identity 
			$(I+tM)\det(I+tM)^{-1}(I-t(M-\rho_MI)=I$ 
			is similarly verified. 
			\qed

		\subsection{The main theorem}
		
		\subsubsection{The main theorem}\label{subsubsection : The main theorem}
		Let $\Delta=(V, {\cal A})$ be a finite digraph which may have multi-arcs and multi-loops, and 
		$\theta^{\rm G}$ a map ${\cal A}\times{\cal A}\rightarrow R$ 
		of the form 
		$
			\theta^{\rm G}(a, a')
			=
			\tau(a')\delta_{\head(a)\tail(a')}-\upsilon(a')\delta_{a'\in S(a)}
		$ 
		where $\tau$ and $\upsilon$ are maps form $\cal A$ to $R$. 
		Recall ${\cal A}(u, v)={\cal A}_{uv}\cup {\cal A}_{vu}$. 
		Fix a total order $\leq$ on $V$. 
		We may assume $u\leq v$ for ${\cal A}(u, v)$ 
		since the definition of ${\cal A}(u, v)$ has the symmetry in $u$ and $v$. 
		Consider the set 
		$$\Phi_{\Delta}=\{(u, v)\in V\times V\mid {\cal A}(u, v)\neq\emptyset\}.
		$$
		Thus, if $(u, v)\in \Phi_{\Delta}$, then one has $u\leq v$. 
		The set $\Phi_{\Delta}$ divides into three subsets 
		$
			\Phi_{\Delta}
			=
			\Phi_{\Delta}^{(1)}
			\sqcup
			\Phi_{\Delta}^{(2)}
			\sqcup
			\Phi_{\Delta}^{(3)},
		$ 
		where 
		$$
		\begin{array}{l}
			\Phi_{\Delta}^{(1)}
			=
			\{
				(u, v)\in{\Phi_{\Delta}}
				\mid 
				u=v
			\}, \\
			\Phi_{\Delta}^{(2)}
			=
			\{
				(u, v)\in{\Phi_{\Delta}}
				\mid 
				u\neq v, 
				\mbox{${\cal A}_{uv}=\emptyset$ or ${\cal A}_{vu}=\emptyset$}
			\}, \\
			\Phi_{\Delta}^{(3)}
			=
			\{
				(u, v)\in{\Phi_{\Delta}}
				\mid 
				u\neq v, 
				\mbox{${\cal A}_{uv}\neq \emptyset$ and ${\cal A}_{vu}\neq\emptyset$}
			\}.
		\end{array}
		$$
		Thus, for $(u, v)\in\Phi_{\Delta}^{\rm (1)}$, 
		${\cal A}(u, v)={\cal A}(u, u)(={\cal A}_{uv})$ is the set of loops with nest $u$. 
		For $(u, v)\in V\times V$, 
		we consider the following polynomial 
		$$
			f_{(u, v)}(t)
			=			
			\left\{
				\begin{array}{ll}
				{\displaystyle
				1+t\sum_{a\in {\cal A}_{uv}}\upsilon(a)
				}, 
				&\mbox{if $u=v$}, \\
				{\displaystyle
				1-t^2
				\left\{\sum_{a\in{\cal A}_{uv}}\upsilon(a)\right\}
				\left\{\sum_{a\in{\cal A}_{vu}}\upsilon(a)\right\}
				},&\mbox{otherwise,}
				\end{array}
			\right.
		$$
		which is invertible in the ring $R[[t]]$ of formal power series. 
		The following lemma follows immediately from the definition.

		\begin{lem}\label{lem : Properties of f}
			For $(u, v)\in V\times V$, we have:
			\begin{itemize}
			\item[{\rm 1)}]
				$f_{(u, v)}(t)=f_{(v, u)}(t)$, 
			\item[{\rm 2)}]
				$f_{(u, v)}(t)=1$ if $(u, v)\notin\Phi_{\Delta}$.
			\end{itemize}
		\end{lem}

		\vspace*{3mm}
		Set $f_{\Delta}(t)=\prod_{(u, v)\in\Phi_{\Delta}}f_{(u, v)}(t)$. 
		Let 
		$$
			A_{\Delta}(\theta^{\rm G})=(a_{ww'}(\theta^{\rm G}))_{w, w'\in V},
			\quad
			D_{\Delta}(\theta^{\rm G})=(d_{ww'}(\theta^{\rm G}))_{w, w'\in V}
		$$ 
		be square matrices of degree $|V|$, 
		whose entries are respectively given by 
		\begin{equation}\label{equation : a and d}
			\begin{array}{l}
			{\displaystyle 
			a_{ww'}(\theta^{\rm G})
			=
			f_{(w, w')}(t)^{-1}a_{ww'},\quad a_{ww'}:=\sum_{a\in{\cal A}_{ww'}}\tau(a),
			}
			\\
			{\displaystyle
			d_{ww'}(\theta^{\rm G})
			=
			\delta_{ww'}
			\sum_{(u, v)\in\Phi_{\Delta}^{(3)}}
			f_{(u, v)}(t)^{-1}(\delta_{wu}d_{uv}+\delta_{wv}d_{vu}),
			\quad
			d_{uv}
			:=
			\sum_{a\in{\cal A}_{uv},\,a'\in{\cal A}_{vu}}
			\tau(a)\upsilon(a').
			}
			\end{array}
		\end{equation}

		\begin{thm}[The Main Theorem]\label{theorem : The main theorem} 
		For a finite digraph $\Delta$, the reciprocal $Z_{\Delta}(t; \theta^{\rm G})^{-1}$ 
		of the generalized weighted zeta function is given by the polynomial 
		\begin{equation}\label{equation : The main theorem}
			f_{\Delta}(t)
			\det\left(I-tA_{\Delta}(\theta^{\rm G})+t^2D_{\Delta}(\theta^{\rm G})\right).
		\end{equation}
		\end{thm}
		{\it Proof.}
		 By Proposition \ref{prop : The three expressions for the generalized weighted zeta function}, 
		 the inverse of the generalized weighted zeta $Z_{\Delta}(t; \theta^{\rm G})$ is given by 
		 $\det(I-tM_{\Delta}(\theta^{\rm G})),$
		 where $M=M_{\Delta}(\theta^{\rm G})=(\theta^{\rm G}(a, a'))_{a, a'\in{\cal A}}$. 
		 Let 
		 $H=(h_{aa'})_{a, a'\in{\cal A}}$ and $J=(j_{aa'})_{a, a'\in{\cal A}}$ 
		 be square matrices with entries 
		 $
		 	h_{aa'}=\tau(a')\delta_{\head(a)\tail(a')}
		$ 
		and 
		$
			j_{aa'}
			=
		 	\upsilon(a')\delta_{a'\in S(a)}.
		$ 
		Obviously we have $M=H-J$. 
		Let $K=(k_{av})_{a\in {\cal A}, v\in V}$ and 
		$L=(l_{ua'})_{u\in V, a'\in {\cal A}}$ denote the matrices 
		with entries 
		$
			k_{av}=\delta_{\head(a)v}
		$ 
		and 
		$
			l_{ua'}=\tau(a')\delta_{u\tail(a')}
		$ 
		respectively. 
		One can easily see that $H=KL$. 
		Thus we have 
		\begin{equation}\label{equation : M by KL a}
			\det(I-tM)=\det\left(I-t(KL-J)\right)=\det\left((I+tJ)-tKL\right).
		\end{equation}

		For each $(u, v)\in\Phi_{\Delta}$, 
		let $J(u, v)$, $K(u, v)$ and $L(u, v)$ denote the submatrices 
		$$
			\begin{array}{ll}
			J(u, v)=(j_{aa'})_{a, a'\in{\cal A}(u, v)}, & j_{aa'}=\upsilon(a')\delta_{a'\in S(a)},\\
			K(u, v)=(k_{aw})_{a\in{\cal A}(u, v), w\in V}, &k_{aw}=\delta_{\head(a)w}, \\
			L(u, v)=(l_{wa'})_{w\in V, a'\in{\cal A}(u, v)}, & l_{wa'}=\tau(a')\delta_{w\tail(a')}.
			\end{array}
		$$
		Thus we have block partition for the matrices $J, K, L$ by 
		$$
			J=(J(u, v))_{(u,v)\in\Phi_{\Delta}}, 
			K=(K(u, v))_{(u,v)\in\Phi_{\Delta}}, 
			L=(L(u, v))_{(u,v)\in\Phi_{\Delta}}.
		$$
		Note that any total order on $\Phi_{\Delta}$ makes $J=(J(u, v))_{(u,v)\in\Phi_{\Delta}}$ 
		a block diagonal matrix. 
		In what follows, 
		we fix such an order on $\Phi_{\Delta}$. 
		Note that any diagonal block $J(u, v)$ ($(u, v)\in\Phi_{\Delta}$) is column constant. 
		Hence, the matrix $I+tJ$ is invertible, since, 
		by Lemma \ref{lem : The inverse of a column constant matrix}, 
		each diagonal block $I+tJ(u, v)$ is invertible for any $(u, v)\in\Phi_{\Delta}$. 
		Therefore, (\ref{equation : M by KL a}) is written by 
		$
			\det(I+tJ)\det(I-t(I+tJ)^{-1}KL), 
		$
		which equals
		\begin{equation}\label{equation : M by KL b}
			\det(I+tJ)\det(I-tL(I+tJ)^{-1}K), 
		\end{equation}
		since, for two matrices $A$, $B$, 
		we have $\det(I-AB)=\det(I-BA)$ if $AB$ and $BA$ are square matrices. 
		The inverse $(I+tJ)^{-1}$ is also block diagonal, i.e, 
		$$
			(I+tJ)^{-1}
			=
			\bigoplus_{(u, v)\in\Phi_{\Delta}}
			(I+tJ(u, v))^{-1}.
		$$
		Hence we have 
		$
		L(I+tJ)^{-1}K
		=
		\sum_{(u, v)\in\Phi_{\Delta}}L(u, v)(I+tJ(u, v))^{-1}K(u, v).
		$
		We compute 
		$$
		T(u, v):=L(u, v)(I+tJ(u, v))^{-1}K(u, v)
		$$ 
		for each $(u, v)\in{\Phi}_{\Delta}^{(i)}$, $i=1,2,3$.

		If $(u, v)\in\Phi_{\Delta}^{(1)}$, 
		then ${\cal A}(u, v)$ consists of all arcs with nest $u(=v)$, 
		say ${\cal A}(u, u)=\{a_1, \dots, a_k\}$. 
		In this case, 
		the square matrix $J(u, u)$ is given by 
		$$
			\left[
				\begin{array}{ccc}
				\upsilon(a_1) & \cdots & \upsilon(a_k)\\
				\vdots &&\vdots\\
				\upsilon(a_1) & \cdots & \upsilon(a_k)
				\end{array}
			\right]
		$$
		which is column constant. 
		Hence, by Lemma \ref{lem : The determinant of a column constant matrix}, 
		it follows that 
		$
			\det(I+tJ(u, u))
			=
			1+t\sum_{i=1}^k\upsilon(a_i)
			=
			f_{uu}(t).
		$ 
		Thus, by Lemma \ref{lem : The inverse of a column constant matrix}, 
		the inverse $(I+tJ(u, u))^{-1}$ is given by 
		$f_{uu}(t)^{-1}(I-t(J(u, u)-({\rm tr}\,J(u, u)) I))$. 
		For each $a\in{\cal A}(u, u)$ and $v\in V$, 
		the $(a, v)$-entry of the product $(J(u, u)-({\rm tr}\, J(u, u)I)K(u, u)$ 
		is given by 
		\begin{eqnarray*}
			\sum_{a'\in {\cal A}(u, u)}
			\left(\upsilon(a')-{\rm tr}\,J(u, u)\delta_{aa'}\right)\delta_{\head(a')v}
			&=&
			\sum_{a'\in {\cal A}(u, u)}\upsilon(a')\delta_{uv}
			-
			\sum_{a'\in {\cal A}(u, u)}{\rm tr}\,J(u, u)\delta_{aa'}\delta_{uv}\\
			&=&
			\sum_{a'\in {\cal A}(u, u)}\upsilon(a')\delta_{uv}
			-
			{\rm tr}\,J(u, u)\delta_{uv}\\
			&=&0.
		\end{eqnarray*}
		Therefore we have 
		\begin{eqnarray*}
			\left(I+tJ(u, u)\right)^{-1}K(u, u)
			&=&
			f_{uu}(t)^{-1}
			\{I-t(J(u, u)-{\rm tr}\,J(u,u)I)\}K(u, u)
			\\
			&=&
			f_{uu}(t)^{-1}
			K(u, u), 
		\end{eqnarray*}
		and it follows that 
		$
			T(u, v)
			=f_{(u, v)}(t)^{-1}L(u, v)K(u,v)
		$ 
		for $(u, v)\in\Phi_{\Delta}^{(1)}$.

		Suppose that $(u, v)\in\Phi_{\Delta}^{(2)}$. 
		Due to symmetry of the definition of ${\cal A}(u, v)$, 
		we may assume ${\cal A}(u, v)={\cal A}_{uv}$. 
		In this case, we have $J(u, v)=O$. 
		Hence $(I+tJ(u, v))^{-1}=I$. 
		Note also that $f_{(u, v)}(t)=1$ for $(u, v)\in\Phi_{\Delta}^{(2)}$, i.e., 
		we have $\det(I+tJ(u, v))=f_{(u, v)}(t)$ for $(u, v)\in\Phi_{\Delta}^{(2)}$. 
		Therefore we have 
		$
			T(u, v)
			=f_{(u, v)}(t)^{-1}L(u, v)K(u, v)
		$ 
		for $(u, v)\in\Phi_{\Delta}^{(2)}$.

		Let $(u, v)\in\Phi_{\Delta}^{(3)}$. 
		Hence ${\cal A}(u, v)={\cal A}_{uv}\sqcup{\cal A}_{vu}$, 
		${\cal A}_{uv}\neq\emptyset$, ${\cal A}_{vu}\neq\emptyset$; 
		say 
		${\cal A}_{uv}=\{a_1, \dots, a_k\}$, ${\cal A}_{vu}=\{a_{k+1}, \dots, a_{k+l}\}$. 
		If we consider the following two column constant matrices 
		$
			J_1=(\upsilon(a'))_{a\in{\cal A}_{uv}, a'\in{\cal A}_{vu}}
		$, 
		$
			J_2=(\upsilon(a'))_{a\in{\cal A}_{vu}, a'\in{\cal A}_{uv}}
		$, then 
		the matrix $I+tJ(u, v)$ is written by 
		$$
			\left[
			\begin{array}{cc}
				I & tJ_1\\
				tJ_2 & I
			\end{array}
			\right]. 
		$$ 
		If we apply Lemma \ref{lem : The inverse of a column constant matrix}, replacing $t$ by $-t^2$, 
		we can see that $I-(tJ_1)(tJ_2)=I-t^2J_1J_2$ is invertible since the matrix $J_1J_2$ is column constant. 
		Therefore, it follows from Corollary \ref{cor : The inverse of the block matrix with Schur complements} that 
		the matrix $I+tJ(u, v)$ is invertible, and 
		the inverse $(I+tJ(u, v))^{-1}$ is given by 
		\begin{equation}\label{equation : The inverse of I+tJ(u, v)}
			\left[
			\begin{array}{cc}
				(I-t^2J_1J_2)^{-1} & -tJ_1(I-t^2J_2J_1)^{-1}\\
				-t(I-t^2J_2J_1)^{-1}J_2 & (I-t^2J_2J_1)^{-1}
			\end{array}
			\right].
		\end{equation}
		It also follows from Corollary \ref{cor : The inverse of the block matrix with Schur complements} that 
		$
			\det(I+tJ(u, v))
			=
			\det(I-t^2J_1J_2)
			=
			\det(I-t^2J_2J_1).
		$ 
		Since the matrix $-t^2J_1J_2$ is column constant, 
		it follows from Lemma \ref{lem : The determinant of a column constant matrix} that 
		$
			\det(I-t^2J_1J_2)
			=
			1+\rho_{-t^2J_1J_2}
			=
			1-t^2\rho_{J_1}\rho_{J_2}, 
		$ 
		and one can easily see from the definition that this equals $f_{(u, v)}(t)$. 
		Therefore we have $\det(I+tJ(u, v))=f_{(u, v)}(t)$ for $(u, v)\in\Phi_{\Delta}^{(3)}$, 
		which also shows that $\det(I-t^2J_1J_2)=\det(I-t^2J_2J_1)=f_{(u, v)}(t)$. 
		By Lemma \ref{lem : The inverse of a column constant matrix}, 
		the $(1, 1)$-block matrix $(I-t^2J_1J_2)^{-1}$ in (\ref{equation : The inverse of I+tJ(u, v)}) is 
		written by 
		$$
			f_{(u, v)}(t)^{-1}
			\{
				I+t^2(J_1J_2-\rho_{J_1J_2}I)
			\}.
		$$ 
		Similarly we have 
		$
			(I-t^2J_2J_1)^{-1}
			=
			f_{(u, v)}(t)^{-1}
			\{
				I+t^2(J_2J_1-\rho_{J_2J_1}I)
			\}, 
		$
		which gives the $(2, 2)$-block in (\ref{equation : The inverse of I+tJ(u, v)}). 
		The remaining blocks contain the products $J_1J_2J_1$ and $J_2J_1J_2$, 
		and it is easily seen that 
		$J_1J_2J_1=\rho_{J_2}\rho_{J_1}J_1$ 
		and $J_2J_1J_2=\rho_{J_1}\rho_{J_2}J_2$. 
		Thus, $(I+tJ(u, v))^{-1}$ equals
		\begin{equation}\label{equation : (I+tJ(u, v))^{-1}}
			f_{(u, v)}(t)^{-1}
			\left\{
				\left[
				\begin{array}{cc}
					I & O\\
					O & I
				\end{array}
				\right]
				-t
				\left[
				\begin{array}{cc}
					O & J_1\\
					J_2 & O
				\end{array}
				\right]
				+t^2
				\left[
				\begin{array}{cc}
					J_1J_2-\rho_{J_1J_2}I & O\\
					O & J_2J_1-\rho_{J_2J_1}I
				\end{array}
				\right]
			\right\}. 
		\end{equation}

		\noindent
		Hence the matrix $T(u, v)$ is a linear combination of 
		the following three matrices 
		$L(u, v)K(u, v)$, $L(u, v)J(u, v)K(u, v)$ and 
		$
			L(u, v)
			\left\{
				P_1\oplus P_2
			\right\}
			K(u, v), 
		$
		where 
		$
			P_1=J_1J_2-\rho_{J_1J_2}I
		$, 
		$
			P_2=J_2J_1-\rho_{J_2J_1}I. 
		$
		For any two vertices $\xi, \eta\in V$, 
		let 
		$
			K_{\xi\eta}=(k_{aw})_{a\in{\cal A}_{\xi\eta}, w\in V}. 
		$ 
		Since the matrix $K(u, v)$ is partitioned by two blocks $K_{uv}$ and $K_{vu}$, 
		the product 
		$
			\left\{
				P_1\oplus P_2
			\right\}
			K(u, v)
		$ is also partitioned by the following two blocks 
		$
			P_1K_{uv}, 
			P_2K_{vu}.
		$ 
		Since $J_1J_2=\rho_{J_1}(\upsilon(a'))_{a, a'\in{\cal A}_{uv}}$, 
		we have 
		$
			J_1J_2K_{uv}
			=
			\rho_{J_1}
			(\upsilon(a'))_{a, a'\in{\cal A}_{uv}}
			(\delta_{\head(a)w})_{a\in{\cal A}_{uv}, w\in V}, 
		$
		the $(a, w)$-entry of which is given by 
		$
			\sum_{b\in{\cal A}_{uv}}
			\rho_{J_1}\upsilon(b)\delta_{\head(b)w}
			=
			\rho_{J_1}\rho_{J_2}\delta_{\head(a)w}
		$ 
		for $a\in{\cal A}_{uv}$ and $w\in V$. 
		Hence we have $J_1J_2K_{uv}=\rho_{J_1J_2}K_{uv}$, 
		which implies $P_1K_{uv}=O$ 
		since $\rho_{J_1}\rho_{J_2}=\rho_{J_1J_2}$ 
		by Lemma \ref{lem : The product of column constant matrices}. 
		In the same manner, one can show that $P_2K_{vu}=O$. 
		Thus we have 
		$
			L(u, v)
			\left\{
				P_1\oplus P_2
			\right\}
			K(u, v)
			=O, 
		$
		and it follows that 
		$
			T(u, v)
			=
			f_{(u, v)}(t)^{-1}
			\left\{
			L(u, v)K(u, v)-tL(u, v)J(u, v)K(u, v)
			\right\}
		$ 
		for $(u, v)\in\Phi_{\Delta}^{(3)}$. 
		Consequently we have 
		\begin{eqnarray*}
			 \lefteqn{L(I+tJ)^{-1}K}\\
			 &&=
			  \sum_{i=1}^3\sum_{(u, v)\in\Phi^{(i)}_{\Delta}}L(u, v)(I+tJ(u, v))^{-1}K(u, v)\\
			 &&=
			 \sum_{(u, v)\in\Phi_{\Delta}}f_{(u, v)}(t)^{-1}L(u, v)K(u, v)
			 -
			  t\sum_{(u, v)\in\Phi^{(3)}_{\Delta}}f_{(u, v)}(t)^{-1}L(u, v)J(u, v)K(u, v).
		\end{eqnarray*}

		For each $(u, v)\in\Phi_{\Delta}$, 
		let 
		$
			A(u, v):=L(u, v)K(u, v), 
		$ 
		and 
		$
			D(u, v):=L(u, v)J(u, v)K(u, v),
		$
		say $A(u, v)=(\alpha_{ww'}(u, v))_{w, w'\in V}$ and $D(u, v)=(\beta_{ww'}(u, v))_{w, w'\in V}$. 
		Accordingly we have
		\begin{equation}
			 L(I+tJ)^{-1}K
			 =
			  \sum_{(u, v)\in\Phi_{\Delta}}f_{(u, v)}(t)^{-1}A(u, v)
			 -
			  t\sum_{(u, v)\in\Phi^{(3)}_{\Delta}}f_{(u, v)}(t)^{-1}D(u, v).
		\end{equation}

		Let $w, w'\in V$. 
		The $(w, w')$-entry $\alpha_{ww'}(u, v)$ of $A(u, v)$ is given by
		$$
			\sum_{a\in{\cal A}(u, v)}\tau(a)\delta_{w\tail(a)}\delta_{\head(a)w'}.
		$$ 
		If $(u, v)=(u, u)\in\Phi_{\Delta}^{(1)}$, 
		then ${\cal A}(u, u)={\cal A}_{uu}$ and we have 
		$
			\alpha_{ww'}(u, v)
			=
			\sum_{a\in{\cal A}_{uu}}
			\tau(a)\delta_{wu}\delta_{uw'}.
		$ 
		This equals 
		$
			\delta_{wu}\delta_{uw'}
			\sum_{a\in{\cal A}_{uu}}
			\tau(a)
			=
			\delta_{wu}\delta_{uw'}a_{uu},
		$ 
		where $a_{uv}$ is defined in (\ref{equation : a and d}). 
		In the case where $(u, v)\in \Phi_{\Delta}^{(2)}$, 
		then ${\cal A}_{uv}=\emptyset$ or ${\cal A}_{vu}=\emptyset$. 
		In this case where ${\cal A}_{vu}=\emptyset$, 
		one has 
		$
			\alpha_{ww'}(u, v)
			=
			\sum_{a\in{\cal A}_{uv}}
			\tau(a)\delta_{wu}\delta_{vw'}, 
		$
		and this equals 
		$
			\delta_{wu}\delta_{w'v}a_{uv}.
		$ 
		Similarly, 
		if ${\cal A}_{uv}=\emptyset$, 
		then one has $\alpha_{ww'}(u, v)=\delta_{wv}\delta_{w'u}a_{vu}.$ 
		For $(u, v)\in\Phi_{\Delta}^{(3)}$, 
		both ${\cal A}_{uv}$ and ${\cal A}_{vu}$ are not empty set, 
		and we have 
		$
			\alpha_{ww'}(u, v)
			=
			\delta_{wu}\delta_{w'v}a_{uv}
			+
			\delta_{wv}\delta_{w'u}a_{vu}.
		$ 
		Putting all these together, 
		one can see that the $(w, w')$-entry of the matrix 
		$
			\sum_{(u, v)\in\Phi_{\Delta}}f_{(u, v)}(t)^{-1}A(u, v)
		$ 
		is given by 
		$
			f_{(w, w')}(t)^{-1}a_{ww'}.
		$ 
		Therefore we have 
		$$
			A_{\Delta}(\theta^{\rm G})= \sum_{(u, v)\in\Phi_{\Delta}}f_{(u, v)}(t)^{-1}A(u, v).
		$$
		On the other hand, for $(u, v)\in\Phi_{\Delta}^{(3)}$, 
		the $(w, w')$-entry $\beta_{ww'}(u, v)$ of $D(u, v)$ is given by 
		$$
			\sum_{a, a'\in{\cal A}(u, v)}
			\tau(a)\delta_{w\tail(a)}\upsilon(a')\delta_{a'\in S(a)}\delta_{\head(a')w'}.
		$$
		One can easily see that this equals 
		$
			\delta_{wu}\delta_{w'u}d_{uv}
			+
			\delta_{wv}\delta_{w'v}d_{vu}.
		$ 
		Thus the $(w, w')$-entry $d_{ww'}(\theta^{\rm G})$ 
		of the matrix $\sum_{(u, v)\in\Phi^{(3)}_{\Delta}}f_{(u, v)}(t)^{-1}D(u, v)$ 
		is 
		$
			\sum_{(u, v)\in\Phi_{\Delta}^{(3)}}
			f_{(u, v)}(t)^{-1}
			(
			\delta_{wu}\delta_{w'u}d_{uv}
			+
			\delta_{wv}\delta_{w'v}d_{vu}
			)
		$, 
		and this equals 
		$
			\delta_{ww'}
			\sum_{(u, v)\in\Phi_{\Delta}^{(3)}}
			f_{(u, v)}(t)^{-1}
			(
				\delta_{wu}d_{uv}+\delta_{wv}d_{vu}
			),
		$ 
		which shows that 
		$$
			 D_{\Delta}(\theta^{\rm G})=\sum_{(u, v)\in\Phi^{(3)}_{\Delta}}f_{(u, v)}(t)^{-1}D(u, v). 
		$$

		Finally, 
		since we have verified that $\det(I+tJ(u, v))=f_{(u, v)}(t)$ for each $(u, v)\in\Phi_{\Delta}$, 
		it follows that $\det(I+tJ)=\prod_{(u, v)\in\Phi_{\Delta}}f_{(u, v)}(t)=f_{\Delta}(t)$. 
		Therefore, from (\ref{equation : M by KL}), 
		we have 
		$$
			\det(I-tM)=f_{\Delta}(t)\det(I-tA_{\Delta}(\theta^{\rm G})+t^2D_{\Delta}(\theta^{\rm G})).
		$$\qed

		\vspace*{3mm}
		The identity for the generalized weighted zeta $Z_{\Delta}(t; \theta^{\rm G})$ 
		verified in Theorem \ref{theorem : The main theorem} is called 
		the {\it Ihara expression} for $Z_{\Delta}(t; \theta^{\rm G})$, 
		and the matrices $A_{\Delta}(\theta^{\rm G})$ and $D_{\Delta}(\theta^{\rm G})$ 
		are called the {\it weighted adjacency matrix} and the {\it weighted backtrack matrix} 
		for $Z_{\Delta}(t; \theta^{\rm G})$, respectively.

		\subsubsection{Comments on the adjacency matrix and the backtrack matrix}
		\label{subsubsection : Comments on the adjacency matrix and the backtrack matrix}

		The weighted adjacency matrix $A_{\Delta}(\theta^{\rm G})$ and 
		the weighted backtrack matrix $D_{\Delta}(\theta^{\rm G})$ are 
		practically discrepant from the ordinary ones. 
		We will see the difference between them for the case of the Ihara zeta $Z_{\Delta(\Gamma)}(t; \theta^{\rm I})$ 
		for a finite simple graph $\Gamma$. 
		Let $\Gamma=(V, E)$ be a finite simple graph, 
		where $V$ is the vertex set, $E$ the edge set. 
		Suppose that a digraph $\Delta=\Delta(\Gamma)$ is the symmetric digraph of $\Gamma$. 
		Let ${\cal A}$ denote the arc set of $\Delta$. 
		$\Phi_{\Delta}^{(i)}$ ($i=1, 2, 3$) is defined as in \ref{subsubsection : The main theorem}.

		\begin{lem}\label{lem : Phi_Delta for finite simple graphs}
		Let $\Delta=(V, {\cal A})$ be the symmetric digraph of a 
		finite simple graph. Then we have:
		{\rm 1)}\,$\Phi_{\Delta}^{(1)}=\Phi_{\Delta}^{(2)}=\emptyset,$
		{\rm 2)}\,${\cal A}(u, v)\neq\emptyset$ implies $|{\cal A}(u, v)|=2$. 
				In this case, we have $|{\cal A}_{uv}|=|{\cal A}_{vu}|=1$. 
		\end{lem}


		It follows from Lemma \ref{lem : Phi_Delta for finite simple graphs} that 
		$f_{(u, v)}(t)=1$ for $(u, v)\in\Phi_{\Delta}^{(i)}$ ($i=1, 2$), 
		and $f_{(u, v)}(t)=1-t^2$ for $(u, v)\in\Phi_{\Delta}^{(3)}$. 
		Remark that the coefficients $f_{(u, v)}(t)$ do not 
		depend on $(u, v)\in\Phi_{\Delta}^{(3)}$ in this case. 
		Hence we have 
		$$
			A_{\Delta}(\theta)
			=
			(1-t^2)^{-1}
			\sum_{(u, v)\in\Phi_{\Delta}^{(3)}}
			A(u, v), 
			\quad
			D_{\Delta}(\theta)
			=
			(1-t^2)^{-1}
			\sum_{(u, v)\in\Phi_{\Delta}^{(3)}}
			D(u, v).
		$$ 
		Let $(u, v)\in\Phi_{\Delta}^{(3)}$. 
		Recall that $\theta^{\rm I}$ is obtained by letting $\tau=\upsilon=1$ for $\theta^{\rm G}$. 
		Thus the $(w, w')$-entry $\alpha_{ww'}(u, v)$ of the matrix $A(u, v)$ is given by 
		$
			\alpha_{ww'}(u, v)
			=
			\sum_{a\in{\cal A}(u, v)}\tau(a)\delta_{w\tail(a)}\delta_{\head(a)w'}
			=
			\sum_{a\in{\cal A}_{uv}}\delta_{wu}\delta_{vw'}
			+
			\sum_{a\in{\cal A}_{vu}}\delta_{wv}\delta_{uw'}
			=
			\delta_{wu}\delta_{w'v}
			|{\cal A}_{uv}|
			+
			\delta_{wv}\delta_{w'u}
			|{\cal A}_{vu}|
			=
			\delta_{wu}\delta_{w'v}
			+
			\delta_{wv}\delta_{w'u}. 
		$ 
		This shows that the matrix 
		$
			\sum_{(u, v)\in\Phi_{\Delta}^{(3)}}
			A(u, v)
		$ 
		is nothing but the usual adjacency matrix $A_{\Gamma}$ for $\Gamma$. 
		It also follows from Lemma \ref{lem : Phi_Delta for finite simple graphs} that 
		$
			\beta_{ww'}(u, v)
			=
			\sum_{a, a'\in{\cal A}(u, v)}
			\delta_{w\tail(a)}\delta_{a'\in S(a)}\delta_{\head(a')w'}
			=
			\sum_{a\in{\cal A}_{uv}, a'{\cal A}_{vu}}
			\delta_{wu}\delta_{uw'}
			+
			\sum_{a\in{\cal A}_{vu}, a'{\cal A}_{uv}}
			\delta_{wv}\delta_{vw'}
			=
			\delta_{w=w'=u}
			+
			\delta_{w=w'=v}. 
		$
		This shows that the matrix 
		$
			\sum_{(u, v)\in\Phi_{\Delta}^{(3)}}
			D(u, v)
		$ 
		gives the usual degree matrix $D_{\Gamma}$ of $\Gamma$. 
		Therefore we have 
		$$
			A_{\Delta}(\theta)
			=
			(1-t^2)^{-1}A_{\Gamma},
			\quad
			D_{\Delta}(\theta)
			=
			(1-t^2)^{-1}D_{\Gamma}. 
		$$
		Note that the cardinality $|\Phi_{\Delta}^{(3)}|$ equals 
		the number $|E|$ of edges of $\Gamma$. 
		We also have $f_{(u, v)}(t)=1-t^2$ for each $(u, v)\in\Phi_{\Delta}^{(3)}$,  
		and it follows that $f_{\Delta}(t)=(1-t^2)^{|E|}$. 
		Since the degrees of the square matrix $A_{\Gamma}$ and $D_{\Gamma}$ equal  
		the number $|V|$ of vertices in $\Gamma$, 
		it follows from (\ref{equation : The main theorem}) that 
		$$
			Z_{\Delta}(t)^{-1}
			=
			(1-t^2)^{|E|-|V|}
			\det(I-tA_{\Gamma}+t^2(D_{\Gamma}-I)),
		$$
		which is the classical Bass-Ihara theorem \cite{bass92, ihara66} 
		for a finite simple graph.

		\subsection{Example}
		
\begin{figure}[tbp]
		\center
		\includegraphics[width=4.5cm]{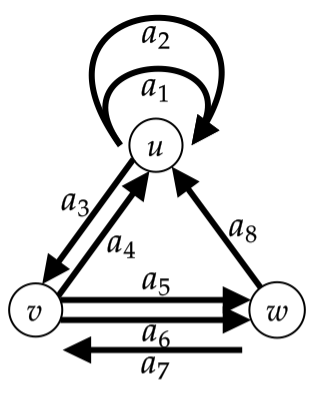}
		\caption{A digraph $\Delta=(V,{\cal A})$}
		\label{figure : graph}
		\end{figure}

		Let $\Delta=(V, {\cal A})$ be a digraph in Figure \ref{figure : graph}, 
		where 
		$
		V=\{ v_1,v_2,v_3 \}
		$ 
		is the vertex set 
		and 
		$
		{\cal A}= \{ a_1, a_2, a_3, a_4, a_5, a_6, a_7, a_8 \}
		$ 
		the arc set. 
		The total order for $V$ is given by $v_1<v_2<v_3$, 
		and for ${\cal A}$ by $a_1<a_2< a_3< a_4<a_5<a_6<a_7<a_8$. 
		Thus we have 
		$$
			\Phi_{\Delta}
			=
			\{
				(v_1, v_1), (v_1, v_2), (v_1, v_3), (v_2, v_3)
			\},
		$$ 
		where 
		$
			\Phi_{\Delta}^{(1)}
			=
			\{
				(v_1, v_1)
			\}, 
			\Phi_{\Delta}^{(2)}
			=
			\{
				(v_1, v_3)
			\}, 
			\Phi_{\Delta}^{(3)}
			=
			\{
				(v_1, v_2), (v_2, v_3)
			\}. 
		$
		The arc sets ${\cal A}_{uv}$'s are given by 
		$
			{\cal A}_{v_1v_1}=\{a_1, a_2\}, 
			{\cal A}_{v_1v_2}=\{a_3\}, 
			{\cal A}_{v_1v_3}=\emptyset, 
			{\cal A}_{v_2v_1}=\{a_4\}, 
			{\cal A}_{v_2v_3}=\{a_5, a_6\}, 
			{\cal A}_{v_3v_2}=\{a_7\}.
		$ 
		The matrices $J,K,L$ and $J(u, v), K(u, v), L(u, v)$ for each $(u, v)\in\Phi_{\Delta}$ are given by

		\begin{eqnarray*}
		&&J =
		\left[
		\begin{array}{cccccccc}
		\upsilon(a_1) & \upsilon(a_2) & 0 & 0 & 0 & 0 & 0 & 0 \\
		\upsilon(a_1) & \upsilon(a_2) & 0 & 0 & 0 & 0 & 0 & 0 \\
		0 & 0 & 0 & \upsilon(a_4) & 0 & 0 & 0 & 0 \\
		0 & 0 & \upsilon(a_3) & 0 & 0 & 0 & 0 & 0 \\
		0 & 0 & 0 & 0 & 0 & 0 & \upsilon(a_7) & 0 \\
		0 & 0 & 0 & 0 & 0 & 0 & \upsilon(a_7) & 0 \\
		0 & 0 & 0 & 0 & \upsilon(a_5) & \upsilon(a_6) & 0 & 0 \\
		0 & 0 & 0 & 0 & 0 & 0 & 0 & 0 
		\end{array}
		\right],
		\quad
		K =
		\left[
		\begin{array}{ccc}
		1 & 0 & 0 \\
		1 & 0 & 0 \\
		0 & 1 & 0 \\
		1 & 0 & 0 \\
		0 & 0 & 1 \\
		0 & 0 & 1 \\
		0 & 1 & 0 \\
		1 & 0 & 0 
		\end{array}
		\right],
		\\
		&&L =
		\left[
		\begin{array}{cccccccc}
		\tau(a_1) & \tau(a_2) & \tau(a_3) & 0 & 0 & 0 & 0 & 0 \\
		0 & 0 & 0 & \tau(a_4) & \tau(a_5) & \tau(a_6) & 0 & 0 \\
		0 & 0 & 0 & 0 & 0 & 0 & \tau(a_7) & \tau(a_8) 
		\end{array}
		\right],
		\\
		&&
		J(v_1,v_1) = \begin{bmatrix}
			\upsilon(a_1) & \upsilon(a_2) \\
			\upsilon(a_1) & \upsilon(a_2) 
		\end{bmatrix},
		\quad
		K(v_1,v_1) = \begin{bmatrix}
			1 & 0 & 0\\
			1 & 0 & 0
		\end{bmatrix},
		\quad
		L(v_1,v_1) = \begin{bmatrix}
			\tau(a_1) & \tau(a_2)\\
			0 & 0\\
			0 & 0
		\end{bmatrix}, 
		\\
		&&
		J(v_1,v_2) = \begin{bmatrix}
			0& \upsilon(a_4) \\
			\upsilon(a_3) & 0
		\end{bmatrix},
		\quad
		K(v_1,v_2) = \begin{bmatrix}
		0 & 1 & 0\\
		1 & 0 & 0
		\end{bmatrix}, 
		\quad
		L(v_1,v_2) = \begin{bmatrix}
			\tau(a_3) & 0\\
			0 & \tau(a_4)\\
			0 & 0
		\end{bmatrix}, 
		\\
		&&
		J(v_2,v_3) = \begin{bmatrix}
		0 & 0 & \upsilon(a_7) \\
		0 & 0 & \upsilon(a_7) \\
		\upsilon(a_5)& \upsilon(a_6) & 0
		\end{bmatrix},
		K(v_2,v_3) = \begin{bmatrix}
			0 & 0 & 1\\
			0 & 0 & 1\\
			0 & 1 & 0\\
		\end{bmatrix},
		L(v_2,v_3) 
			=
			 \begin{bmatrix}
				0 & 0 & 0\\
				\tau(a_5) & \tau(a_6) & 0\\
				0 & 0 & \upsilon(a_7)
			\end{bmatrix}, 
		\\
		&&
			J(v_1,v_3) 
			=
			 \begin{bmatrix}
				0
			\end{bmatrix},
			\quad
			K(v_1,v_3) 
			=
			 \begin{bmatrix}
				1 & 0 & 0\\
			\end{bmatrix}, 
			\quad
			L(v_1,v_3) 
			=
			\begin{bmatrix}
				0\\
				0\\
				\tau(a_8) 
			\end{bmatrix}.
		\end{eqnarray*}
		The matrices $A(u, v)$ and $D(u, v)$ are for example 
		\begin{eqnarray*}
		&&
			A(v_1, v_1)
			=
			\left[
				\begin{array}{ccc}
					\tau(a_1)+\tau(a_2) & 0 & 0\\
					0 & 0 & 0\\
					0 & 0 & 0
				\end{array}
			\right],
			\quad
			A(v_1, v_2)
			=
			\left[
				\begin{array}{ccc}
					0 & \tau(a_3) & 0\\
					\tau(a_4) & 0 & 0\\
					0 & 0 & 0
				\end{array}
			\right],
		\\
		&&
			A(v_1, v_3)
			=
			\left[
				\begin{array}{ccc}
					0 & 0 & 0\\
					0 & 0 & 0\\
					\tau(a_8) & 0 & 0
				\end{array}
			\right],
		\\
		&&
			D(v_1, v_2)
			=
			\left[
				\begin{array}{ccc}
					\tau(a_3)\upsilon(a_4) & 0 & 0\\
					0 & \tau(a_4)\upsilon(a_3) & 0\\
					0 & 0 & 0
				\end{array}
			\right],
		\\
		&&
			D(v_2, v_3)
			=
			\left[
				\begin{array}{ccc}
					0 & 0 & 0\\
					0 & (\tau(a_5)+\tau(a_6))\upsilon(a_7) & 0\\
					0 & 0 & \tau(a_7)(\upsilon(a_5)+\upsilon(a_6))
				\end{array}
			\right],
		\end{eqnarray*}
		etc. 
		The polynomials $f_{(u, v)}(t)$ are given by 
		$
			f_{(v_1, v_1)}(t)
			=
			1+t(\upsilon(a_1)+\upsilon(a_2)),
			f_{(v_1, v_2)}(t)
			=
			1-t^2\upsilon(a_3)\upsilon(a_4),
			f_{(v_1, v_3)}(t)
			=
			1, 
			f_{(v_2, v_3)}(t)
			=
			1-t^2(\upsilon(a_5)+\upsilon(a_6))\upsilon(a_7), 
		$
		and we have 
		\begin{eqnarray*}
			A_{\Delta}(\theta^{\rm G}) &=& 
				\left[
				\begin{array}{ccc}
				\frac{\tau(a_1)+\tau(a_2)}{1+t(\upsilon(a_1)+\upsilon(a_2)} & \frac{\tau(a_3)}{1-t^2\upsilon(a_3)\upsilon(a_4)} & 0 \\
				\frac{\tau(a_4)}{1-t^2\upsilon(a_3)\upsilon(a_4)} & 0 & \frac{\tau(a_5)+\tau(a_6)}{1-t^2(\upsilon(a_5)+\upsilon(a_6))							\upsilon(a_7)}  \\
				\tau(a_8) & \frac{\tau(a_7)}{1-t^2(\upsilon(a_5)+\upsilon(a_6))\upsilon(a_7)}  & 0
			\end{array}
			\right]
			\\
			D_{\Delta}(\theta^{\rm G}) &=& 
				\left[
				\begin{array}{ccc}
				 \frac{\tau(a_3)\upsilon(a_4)}{1-t^2\upsilon(a_3)\upsilon(a_4)}  & 0 & 0 \\
				0 &  \frac{\tau(a_4)\upsilon(a_3)}{1-t^2\upsilon(a_3)\upsilon(a_4)}  + \frac{(\tau(a_5)+\tau(a_6))\upsilon(a_7)}{1-								t^2(\upsilon(a_5)+\upsilon(a_6))\upsilon(a_7)} & 0\\
				0 & 0 &  \frac{\tau(a_7)(\upsilon(a_5)+\upsilon(a_6))}{1-t^2(\upsilon(a_5)+\upsilon(a_6))\upsilon(a_7)} 
			\end{array}
			\right].
		\end{eqnarray*}
		The main theorem shows that 
		$
			\det(I-tM_{\Delta}(\theta^{\rm G}))
			=
			f_{\Delta}(t)
			\det(I-tA_{\Delta}(\theta^{\rm G})+t^2D_{\Delta}(\theta^{\rm G})),
		$ 
		where 
		$
			f_{\Delta}(t)
			=
			\prod_{(u, v)\in\Phi_{\Delta}}f_{(u, v)}(t)
		$ 
		and 
		$
			M_{\Delta}(\theta^{\rm G})
			=
			(\theta^{\rm G}(a, a'))_{a, a'\in{\cal A}}. 
		$

			\vspace*{1cm}
			\noindent{\em Acknowledgement.}
			\vspace{0.2cm}

			A. Ishikawa is partially supported by Grant-in-Aid for JSPS Fellows (Grant No. 20J20590).

\end{document}